\documentclass[11pt]{article}

\usepackage[margin= 2.2 cm]{geometry}
\usepackage[utf8]{inputenc} 
\usepackage{amsmath}
\usepackage{graphicx}
\usepackage{subcaption}
\usepackage{amsfonts}
\usepackage{mathtools}
\usepackage{subcaption}
\usepackage{amssymb}
\usepackage{amsthm}
\usepackage{comment}
\usepackage{chngpage}
\usepackage{hyperref}
\hypersetup{
	colorlinks=true,
	linkcolor=black,
	citecolor=black,
	filecolor=black,      
	urlcolor=black,
}
\usepackage{enumerate}
\usepackage{makeidx}
\usepackage{thmtools}
\usepackage{thm-restate}

\usepackage{tikz}
\usetikzlibrary{positioning,chains,fit,shapes,calc,patterns}

\makeindex

\makeatletter
\renewenvironment{proof}[1][\proofname] {\par\pushQED{\qed}\normalfont\topsep6\p@\@plus6\p@\relax\trivlist\item[\hskip\labelsep\bfseries#1\@addpunct{.}]\ignorespaces}{\popQED\endtrivlist\@endpefalse}
\makeatother

\newtheorem{proposition}{Proposition}[section]

\newtheorem{lemma}[proposition]{Lemma}
\newtheorem{theorem}[proposition]{Theorem}

\theoremstyle{definition}
\newtheorem{definition}[proposition]{Definition}

\newtheorem*{remark*}{Remark}
\newtheorem*{theorem*}{Theorem}
\newtheorem{claim}[proposition]{Claim}

\usepackage[shortlabels]{enumitem}
\usepackage[capitalise]{cleveref}
\usepackage{url}

\usepackage[numbers,sort]{natbib}

\renewcommand{\AA}{\mathcal{A}}
\newcommand{\MM}{\mathcal{M}}
\newcommand{\PP}{\mathbb{P}}

\newcommand{\ex}{\operatorname{ex}}
\newcommand{\XX}{\mathbf{X}}
\newcommand{\YY}{\mathbf{Y}}
\newcommand{\HH}{\mathcal{H}}

\newcommand{\eps}{\varepsilon}
\newcommand{\GG}{\mathcal{G}}
\newcommand{\vu}{\mathbf{u}}
\newcommand{\vv}{\mathbf{v}}
\newcommand{\Hom}{\mathrm{Hom}}

\title{Small subgraphs of large minimum degree}
\author{Zach Hunter\footnote{ETH Zurich, Switzerland, {\tt zach.hunter@math.ethz.ch}. Research supported in part by SNSF grant 200021-228014.} \and Barnab\'as Janzer\footnote{Mathematical Institute, University of Oxford, UK, {\tt barnabas.janzer@magd.ox.ac.uk}} \and Oliver Janzer\footnote{Institute of Mathematics, EPFL, Lausanne, Switzerland, {\tt oliver.janzer@epfl.ch}} \and Benny Sudakov\footnote{ETH Zurich, Switzerland, {\tt benjamin.sudakov@math.ethz.ch}. Research supported in part by SNSF grant 200021-228014.}}

\date{}

\begin{document}

\maketitle

\begin{abstract}
    In 1991, Erd\H os, Faudree, Rousseau and Schelp asked to estimate the smallest $u=u(n,d,s)$ such that every $n$-vertex graph with average degree at least $d$ has a subgraph on at most $u$ vertices with minimum degree at least $s$. We resolve this problem up to polylogarithmic factors for all even values of $s$. More precisely, we prove that every $n$-vertex graph with minimum degree $d\geq s$ has a subgraph on $O_s\!\left(\max\left\{\frac{n}{d^{s/(s-2)}}(\log n)^{180},\log n\right\}\right)$ vertices with minimum degree at least $s$, whenever $s\ge 4$ is even. This improves a result of Janzer, Sudakov and Tomon, who obtained a similar conclusion for subgraphs of \emph{average} degree at least $s$.

    We also prove that for every positive integer $s$ and every $\eps>0$, there exists a finite family $\mathcal{H}$ of graphs with minimum degree $s$ such that $\ex(n,\mathcal{H})=O(n^{2-2/s+\eps})$. This is tight and is the first result of similar kind that applies to odd values of $s$ greater than 3.
\end{abstract}

\section{Introduction}

A basic theme in extremal graph theory is to understand how global density forces local structure.  One of the simplest manifestations of this principle is the familiar fact that every graph of average degree at least $2s$ contains a non-empty subgraph of minimum degree at least $s$: one simply removes vertices of degree less than $s$ until this is no longer possible.  This argument, however, gives no useful control on the number of vertices that remain.  In this paper we study the corresponding quantitative problem: how small can one force a subgraph whose \emph{every} vertex has large degree?

This question was raised by Erd\H{o}s, Faudree, Rousseau and Schelp in a sequence of papers on subgraphs of prescribed minimum degree~\cite{erdos1990subgraphs,erdos1991edge}; see also~\cite{sauermann2019proof} for later developments.  Given positive integers $n$ and $s$, and a real number $d\ge 2s$, let $u(n,d,s)$ denote the smallest integer $u$ such that every $n$-vertex graph of average degree at least $d$ contains a subgraph on at most $u$ vertices with minimum degree at least $s$. By the above discussion
such $u$ always exists. 

The case $s=2$ is exactly the classical girth problem.  Indeed, for every graph $G$, the girth of $G$ is precisely the smallest order of a non-empty subgraph of $G$ with minimum degree at least $2$: a shortest cycle is such a subgraph, while every subgraph of minimum degree at least $2$ contains a cycle.  It is therefore natural to view the problem above as a higher-degree analogue of girth, or generalized girth in the terminology of Feige and Wagner~\cite{feige2016generalized}.  Related generalized-girth questions also arise for hypergraphs.  In particular, Feige's small even-cover conjecture~\cite{feige2008small} predicted that an $r$-uniform $n$-vertex hypergraph with $\Theta(\beta n)$ edges contains an even cover of size $O(n\log n/\beta^{2/(r-2)})$. Guruswami, Kothari and Manohar~\cite{guruswami2022algorithms} proved the conjecture up to polylogarithmic factors; see also~\cite{hsieh2023simple,hsieh2025small} for simpler and sharper proofs. Very recently, Feige's conjecture was resolved in full, independently by Bandeira, Kunisky, Nizi\'c-Nikolac, Pesenti and Wang~\cite{bandeira2026hypergraph} and by Schmidhuber and Hastings~\cite{schmidhuber2026spectral}. 

A closely related problem asks for a small subgraph of large \emph{average} degree.  Let $t(n,d,s)$ be the smallest $t$ such that every $n$-vertex graph of average degree at least $d$ contains a subgraph on at most $t$ vertices with average degree at least $s$.  Feige and Wagner~\cite{feige2016generalized} conjectured that for every fixed $s>2$ and $s\leq d\leq n^{1-2/s}$, we have $t(n,d,s)\leq n d^{-s/(s-2)}(\log n)^{O_s(1)}$.  This was proved by Janzer, Sudakov and Tomon~\cite{janzer2024small}: in fact they established the slightly stronger result that for $s\le d\le n^{1-2/s}$,
every $n$-vertex graph of average degree at least $d$ contains a subgraph of average degree at least $s$ on at most $n d^{-s/(s-2)}(\log d)^{O_s(1)}$
vertices.  A random graph construction gives the matching lower bound $t(n,d,s)=\Omega_s(nd^{-s/(s-2)})$, and hence this determines $t(n,d,s)$ up to polylogarithmic factors. Very recently, the result of Janzer, Sudakov and Tomon was extended to hypergraphs by Kuang and Wang~\cite{kuang2026phase} by establishing a reduction to the graph case.

The minimum degree problem studied in this paper is substantially more difficult. The elementary pruning argument which passes from large average degree to large minimum degree loses a factor of two in the target degree.  Combining that argument with the average degree result above, now applied with target average degree $2s$, gives 
only
$u(n,d,s)\le n d^{-s/(s-1)}(\log d)^{O_s(1)}$, a substantially weaker bound than what we have for $t(n,d,s)$.
On the other hand, we do not know of any better lower bound than in the average degree problem, so one may expect the answer to be the same, that is,
$u(n,d,s)\leq n d^{-s/(s-2)}(\log n)^{O_s(1)}$. Indeed, this bound was conjectured in \cite{janzer2024small}. 

Our first main result proves this conjecture for every even target degree. 

\begin{theorem}\label{thm:even min degree}
Whenever $s\ge 4$ is even and $G$ is an $n$-vertex graph with minimum degree at least $d\ge s$, then $G$ has a subgraph on
\[
O_s\!\left(\max\left\{\frac{n}{d^{s/(s-2)}}(\log n)^{180},\log n\right\}\right)
\]
vertices with minimum degree at least $s$.
\end{theorem}

Random graphs show that \cref{thm:even min degree} is optimal up to the polylogarithmic factor. Moreover, the deletion method shows that a polylogarithmic factor is necessary around the critical scale $d\approx n^{1-2/s}$. On the other hand, when $d$ is polynomially larger than $n^{1-2/s}$, the deletion method does not even rule out the existence of a suitable subgraph of constant size.
When the order of the desired subgraph is required to be constant independent of $n$, the question becomes a classical Tur\'an-type problem. Recall that for a family of graphs $\mathcal H$, the Tur\'an number $\ex(n,\mathcal H)$ is the maximum number of edges in an $n$-vertex graph containing no member of $\mathcal H$. Thus, asking whether average degree at least $n^{1-2/s+\eps}$ forces a bounded-order subgraph of minimum degree at least $s$ is equivalent to asking whether there exists a finite family $\mathcal H$ of graphs of minimum degree at least $s$ such that $\ex(n,\mathcal H)=O\left(n^{2-2/s+\eps}\right)$.

In this direction, Verstra\"ete conjectured that for every $s>2$ and every $\eps>0$, average degree at least $n^{1-2/s+\eps}$ forces a bounded-order subgraph of \emph{average degree} at least $s$.  Jiang and Newman~\cite{jiang2017small} proved this when $s$ is an integer, and Janzer, Sudakov and Tomon~\cite{janzer2024small} subsequently proved it for every real $s>2$.
Jiang and Newman also obtained the corresponding bounded-order \emph{minimum degree} result when $s$ is even.  Janzer \cite{janzer2023rainbow,janzer2023disproof} later proved the stronger result that for any even $s$ or $s=3$ and any $\eps>0$, there is an $s$-regular graph $H$ with $\ex(n,H)=O(n^{2-2/s+\eps})$. Thus, before the present work, the bounded-order minimum-degree problem was known for every even $s$ and for $s=3$, but remained open for every odd $s\ge5$.

Our second main result resolves the remaining cases in the bounded-order regime.

\begin{theorem}\label{thm:small odd subgraph}
Let $s\ge5$ be odd.  Then for every $\eps>0$ there exist $t=t(s,\eps)$ and $n_0=n_0(s,\eps)$ such that, for every $n\ge n_0$, every $n$-vertex graph with average degree greater than $n^{1-2/s+\eps}$ contains a subgraph of minimum degree at least $s$ on at most $t$ vertices.
\end{theorem}

Equivalently, we obtain the following formulation in terms of Tur\'an numbers.

\begin{theorem}\label{thm:odd extremal}
Let $s\ge5$ be odd.  Then, for every $\eps>0$, there exists a finite family $\HH$ of graphs with minimum degree at least $s$ such that
\[
        \ex(n,\HH)=O\!\left(n^{2-2/s+\eps}\right).
\]
\end{theorem}

\subsection*{Proof overview}

Although the constructions used for even and odd $s$ are rather different, the two proofs share a common feature. We develop a mechanism, presented in Section \ref{sec:most degrees large}, that allows us to reduce the problem to the easier one of finding a small subgraph in which instead of all vertices having degree at least $s$, we allow a small fraction to have degree less than $s$. The reduction goes as follows. Assume that we have proved the existence of at least one small subgraph with few vertices of degree less than $s$. Applying this result to random vertex samples proves the existence of \emph{many} such subgraphs of the same order. We then use an argument based on sunflowers that combines several of them to remove the exceptional low-degree vertices altogether. Let us explain this argument. Our task is to show that if we have a large collection of induced subgraphs of order $t$, each containing few vertices of degree less than $s$, then we can find a subgraph of minimum degree $s$ on not much more than $t$ vertices. The basic idea is that if we take several copies of these subgraphs which have the same set of exceptional low-degree vertices, then we may hope that all vertices have degree at least $s$ in the union, as the exceptional vertices send edges to each of the subgraphs. However, they may receive the same set of edges in many of these subgraphs, so we have to use a more refined argument.  Using the sunflower lemma of Alweiss, Lovett, Wu and Zhang~\cite{alweiss2021improved}, we find $s+1$ such subgraphs forming a sunflower with the set of low-degree vertices contained in the kernel. In addition, by a slightly unusual application of the sunflower lemma, we can guarantee that for every vertex in the kernel its degrees in the $s+1$ subgraphs are either all equal or all distinct.  If a vertex in the kernel still has degree less than $s$ in the union, its degrees cannot all be distinct; in the equal-degree case the sunflower property implies that all of its neighbours lie in the kernel.  Such a vertex can therefore be deleted simultaneously from every member of the sunflower without affecting the vertices outside the kernel. After this removal, we iterate the process with the remaining sunflower. It is not hard to see that eventually this procedure yields a non-empty subgraph of minimum degree at least $s$. 
With this machinery at our disposal, our remaining task is to find a small subgraph with almost all vertices having degree at least $s$. How we construct this subgraph is different for our two main results.

\begin{figure}[b]
\centering
\begin{tikzpicture}[thin,
	every node/.style={draw,circle},
    every fit/.style={ellipse,draw,text width=1cm}	]

	% the vertices of X_{1}
	\begin{scope}[xshift=0cm,yshift=0cm,start chain=going below,node distance=1.5mm, inner sep=1.3pt]
		\foreach \i in {0,1,2,3}
		\node[%fsnode,
		on chain] (x1\i) {};
	\end{scope}
	% the vertices of X_{2}
	\begin{scope}[xshift=0cm,yshift=2.5cm,start chain=going below,node distance=1.5mm, inner sep=1.3pt]
		\foreach \i in {0,1,2,3}
		\node[%fsnode,
		on chain] (x2\i) {};
	\end{scope}
    
    % the vertices of Y_{1}
	\begin{scope}[xshift=3.5cm,yshift=0cm,start chain=going below,node distance=1.5mm, inner sep=1.3pt]
		\foreach \i in {0,1,2,3}
		\node[%fsnode,
		on chain] (y1\i) {};
	\end{scope}
	% the vertices of Y_{2}
	\begin{scope}[xshift=3.5cm,yshift=2.5cm,start chain=going below,node distance=1.5mm, inner sep=1.3pt]
		\foreach \i in {0,1,2,3}
		\node[%fsnode,
		on chain] (y2\i) {};
	\end{scope}

	% the vertices of Z_{1}
	\begin{scope}[xshift=7cm,yshift=0cm,start chain=going below,node distance=1.5mm, inner sep=1.3pt]
		\foreach \i in {0,1,2,3}
		\node[%fsnode,
		on chain] (z1\i) {};
	\end{scope}
	% the vertices of Z_{2}
	\begin{scope}[xshift=7cm,yshift=2.5cm,start chain=going below,node distance=1.5mm, inner sep=1.3pt]
		\foreach \i in {0,1,2,3}
		\node[%fsnode,
		on chain] (z2\i) {};
	\end{scope}

	% the vertices of U_{1}
	\begin{scope}[xshift=10.5cm,yshift=0cm,start chain=going below,node distance=1.5mm, inner sep=1.3pt]
		\foreach \i in {0,1,2,3}
		\node[%fsnode,
		on chain] (u1\i) {};
	\end{scope}
	% the vertices of U_{2}
	\begin{scope}[xshift=10.5cm,yshift=2.5cm,start chain=going below,node distance=1.5mm, inner sep=1.3pt]
		\foreach \i in {0,1,2,3}
		\node[%fsnode,
		on chain] (u2\i) {};
	\end{scope}

	%labels shown here are not consistent with tikz names
\coordinate (c) at ($(x10)!0.5!(x13)$);
\draw (c) ellipse [x radius=4mm, y radius=8mm];
\node[left=3.1mm, draw=none] at (c) {};
\coordinate (c) at ($(x20)!0.5!(x23)$);
\draw (c) ellipse [x radius=4mm, y radius=8mm];
\node[left=3.1mm, draw=none] at (c) {};

\coordinate (c) at ($(y10)!0.5!(y13)$);
\draw (c) ellipse [x radius=4mm, y radius=8mm];
\node[right=3.1mm, draw=none] at (c) {};
\coordinate (c) at ($(y20)!0.5!(y23)$);
\draw (c) ellipse [x radius=4mm, y radius=8mm];
\node[right=3.1mm, draw=none] at (c) {};

\coordinate (c) at ($(z10)!0.5!(z13)$);
\draw (c) ellipse [x radius=4mm, y radius=8mm];
\node[right=3.1mm, draw=none] at (c) {};
\coordinate (c) at ($(z20)!0.5!(z23)$);
\draw (c) ellipse [x radius=4mm, y radius=8mm];
\node[right=3.1mm, draw=none] at (c) {};

\coordinate (c) at ($(u10)!0.5!(u13)$);
\draw (c) ellipse [x radius=4mm, y radius=8mm];
\node[right=3.1mm, draw=none] at (c) {};
\coordinate (c) at ($(u20)!0.5!(u23)$);
\draw (c) ellipse [x radius=4mm, y radius=8mm];
\node[right=3.1mm, draw=none] at (c) {};

\def\colors{{"blue","red","cyan","black","brown","gray","black"}}

    \pgfmathparse{\colors[1]}
    \let\edgecolor\pgfmathresult
    \draw[\edgecolor] (x10) to (y12);
    \draw[\edgecolor] (x11) to (y13);
    \draw[\edgecolor] (x12) to (y10);
    \draw[\edgecolor] (x13) to (y11);    \pgfmathparse{\colors[0]}
    \let\edgecolor\pgfmathresult
    \draw[\edgecolor] (x20) to (y22);
    \draw[\edgecolor] (x21) to (y20);
    \draw[\edgecolor] (x22) to (y23);
    \draw[\edgecolor] (x23) to (y21);    \pgfmathparse{\colors[2]}
    \let\edgecolor\pgfmathresult
    \draw[\edgecolor] (x10) to (y21);
    \draw[\edgecolor] (x11) to (y23);
    \draw[\edgecolor] (x12) to (y20);
    \draw[\edgecolor] (x13) to (y22);    

    \pgfmathparse{\colors[3]}
    \let\edgecolor\pgfmathresult
    \draw[\edgecolor] (x20) to (y13);
    \draw[\edgecolor] (x21) to (y10);
    \draw[\edgecolor] (x22) to (y11);
    \draw[\edgecolor] (x23) to (y12);    

    \pgfmathparse{\colors[1]}
    \let\edgecolor\pgfmathresult
    \draw[\edgecolor] (y10) to (z10);
    \draw[\edgecolor] (y11) to (z11);
    \draw[\edgecolor] (y12) to (z13);
    \draw[\edgecolor] (y13) to (z12);    \pgfmathparse{\colors[0]}
    \let\edgecolor\pgfmathresult
    \draw[\edgecolor] (y20) to (z23);
    \draw[\edgecolor] (y21) to (z21);
    \draw[\edgecolor] (y22) to (z20);
    \draw[\edgecolor] (y23) to (z22);    \pgfmathparse{\colors[2]}
    \let\edgecolor\pgfmathresult
    \draw[\edgecolor] (y10) to (z20);
    \draw[\edgecolor] (y11) to (z22);
    \draw[\edgecolor] (y12) to (z21);
    \draw[\edgecolor] (y13) to (z23);    

    \pgfmathparse{\colors[3]}
    \let\edgecolor\pgfmathresult
    \draw[\edgecolor] (y20) to (z12);
    \draw[\edgecolor] (y21) to (z11);
    \draw[\edgecolor] (y22) to (z13);
    \draw[\edgecolor] (y23) to (z10);    

    \pgfmathparse{\colors[1]}
    \let\edgecolor\pgfmathresult
    \draw[\edgecolor] (z10) to (u11);
    \draw[\edgecolor] (z11) to (u12);
    \draw[\edgecolor] (z12) to (u13);
    \draw[\edgecolor] (z13) to (u10);    \pgfmathparse{\colors[0]}
    \let\edgecolor\pgfmathresult
    \draw[\edgecolor] (z20) to (u22);
    \draw[\edgecolor] (z21) to (u21);
    \draw[\edgecolor] (z22) to (u20);
    \draw[\edgecolor] (z23) to (u23);    \pgfmathparse{\colors[2]}
    \let\edgecolor\pgfmathresult
    \draw[\edgecolor] (z10) to (u22);
    \draw[\edgecolor] (z11) to (u21);
    \draw[\edgecolor] (z12) to (u20);
    \draw[\edgecolor] (z13) to (u23);    

    \pgfmathparse{\colors[3]}
    \let\edgecolor\pgfmathresult
    \draw[\edgecolor] (z20) to (u11);
    \draw[\edgecolor] (z21) to (u10);
    \draw[\edgecolor] (z22) to (u12);
    \draw[\edgecolor] (z23) to (u13);

    \node[draw=none] at (-1.25,0.75) {$\cdots$};
    \node[draw=none] at (11.75,0.75) {$\cdots$};

\end{tikzpicture}
\caption{The type of graphs used in the even degree case (in the ideal situation), obtained by taking the graph $H_0$ and replacing the edges by arbitrary perfect matchings between $t$-element sets (here, $k=2$ and $t=4$).}    
\label{fig: even case}
\end{figure}
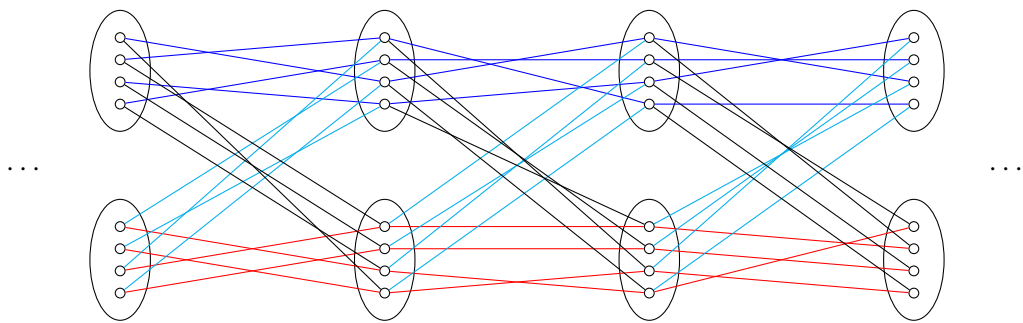

We first describe the even case.  Write $s=2k$. A simple example of a graph with almost all vertices having degree at least $2k$ is the $k$-blowup of a path of large constant length. Let us call this graph $H_0$. This graph is not guaranteed to exist in a very sparse graph as it contains $K_{k,k}$ as a subgraph. Instead, we find a ``matching blowup'' variant of $H_0$ which has significantly more vertices and retains the key property that almost all of its vertices have degree at least $2k$. More precisely, we replace each vertex of $H_0$ by a set of size $t$, and we place a perfect matching (arbitrarily) between any two blobs which were adjacent in $H_0$. See Figure~\ref{fig: even case}. In the ideal picture, the vertex classes we replaced the vertices of $H_0$ with are pairwise disjoint.  Then a typical vertex $v$ receives one distinct neighbour from each of the $2k$ sets that replaced the neighbours of the vertex of $H_0$ whose blob $v$ lives in.  In the actual construction the classes may overlap, so repeated neighbours can lower the degree of some vertices. Details on how we deal with these collisions will be provided at the beginning of Section \ref{sec:even}.

Our main result for odd $s=2k+1$ has a completely different proof. Here we show that any sufficiently dense graph contains one of three possible subgraphs, each of which has few vertices of degree less than $s$. The proof relies on a variant of the homomorphism-counting approach, pioneered by Janzer~\cite{janzer2023rainbow} and developed further by Janzer and Sudakov~\cite{janzer2024hypercube}. The general idea behind this approach is to find inequalities relating the homomorphism counts of certain small subgraphs in our host graph. In a typical application of the homomorphism counting method, one often restricts the count to homomorphisms with extra properties. For example, in coloured problems, in which context the method was first used \cite{janzer2023rainbow}, these may be requirements on the colours of the edges involved. In this paper, we develop this method further and apply it to count homomorphisms satisfying partial injectivity or bounds on the codegrees of the vertices involved. The subgraph we find then depends on the counts of these special homomorphisms. In the ``random-like'' case, when the homomorphism counts do not differ substantially from those in a random graph, we find the graph $P_{\ell}[k]\Box K_2$ as a subgraph, where $P_{\ell}[k]$ denotes the $k$-blowup of the path of length $\ell$ (for some large constant $\ell$) and $\Box$ stands for the Cartesian product. As discussed above, almost all vertices in $P_{\ell}[k]$ have degree $2k$, and after taking the Cartesian product with $K_2$, each vertex receives an extra neighbour. In case the homomorphism counts are not random-like, we exploit the extra structure this guarantees to find one of two different graphs with few vertices of degree less than $s$. We remark that each of those graphs contains $K_{k+1,k+1}$ as a subgraph, so we could not guarantee finding them without the extra information about homomorphism counts, as the edge count of our graph is below the (conjectured) Tur\'an number of $K_{k+1,k+1}$. A more detailed overview of our proof strategy, including the description of the other two subgraphs we might find and the corresponding embedding strategies, will be given at the start of Section \ref{sec:odd}.

\medskip

The rest of the paper is organized as follows.  In Section~\ref{sec:most degrees large} we develop the common reduction from genuine minimum degree to approximate minimum degree, including the sampling and sunflower arguments.  In Section~\ref{sec:even} we prove our result for approximate minimum degree for even $s$ by constructing the matching blow-up described above.  In Section~\ref{sec:odd} we treat odd $s$ and prove the approximate version of Theorem \ref{thm:small odd subgraph}.

\section{Reduction to finding a small subgraph with few vertices of low degree} \label{sec:most degrees large}

In this section, we present the machinery that reduces our task to finding a small subgraph in which \emph{most} (rather than \emph{all}) vertices have degree at least $s$.

Once this reduction is established, the rest of the paper will be devoted to proving the following ``approximate'' versions of Theorems \ref{thm:even min degree} and \ref{thm:small odd subgraph}.
\begin{restatable}{theorem}{restatemuproportion}
\label{thm: mu proportion}
	Let $s\geq 4$ be even, $n$ a positive integer, and let $G$ be an $n$-vertex graph with minimum degree at least $d\geq s$. Then $G$ has a (non-empty) subgraph on $O_{s}\left(\max\left\{\frac{n}{d^{s/(s-2)}}(\log n)^{60},\log n\right\}\right)$ vertices in which the proportion of vertices of degree less than $s$ is at most $1/(\log n)^2$.
\end{restatable}

\begin{theorem} \label{thm:almost all degree s unpruned}
    Let $s\geq 5$ be odd. Then for any $\eps,\mu>0$, there exist some $t=t(s,\mu)$ and $n_0=n_0(s,\eps,\mu)$ such that for each $n\geq n_0$, every $n$-vertex graph with average degree at least $n^{1-2/s+\eps}$ contains a (non-empty) subgraph on at most $t$ vertices in which the proportion of vertices of degree less than $s$ is at most $\mu$.
\end{theorem}

We will prove Theorem \ref{thm: mu proportion} in Section \ref{sec:even} and we will prove Theorem \ref{thm:almost all degree s unpruned} in Section \ref{sec:odd}. The rest of this section is devoted to deducing Theorems \ref{thm:even min degree} and \ref{thm:small odd subgraph} from these results.

\subsection{Many subgraphs where almost all vertices have high degree}

In this section, we use a simple sampling argument to bootstrap Theorems \ref{thm: mu proportion} and \ref{thm:almost all degree s unpruned} to find not just one subgraph where almost every vertex has degree at least $s$, but many such subgraphs of equal orders.

We begin with the counting version of Theorem \ref{thm: mu proportion}.

\begin{lemma}\label{lem:many almost min degree s}
	Let $s\geq 4$ be an even integer, $n$ a positive integer, and let $G$ be an $n$-vertex graph with minimum degree at least $d\geq s$. Then there is some $t=O_{s}(\max\{n(\log n)^{180}/d^{\frac{s}{s-2}},\log n\})$ with $t>0$ such that either $G$ has a subgraph on $t$ vertices with minimum degree at least $s$, or $G$ has at least $(\log n)^{15t}$ induced subgraphs on $t$ vertices in which the proportion of vertices of degree less than $s$ is at most $1/\log n$.
\end{lemma}

\begin{proof}
    We may assume that $n$ is sufficiently large in terms of $s$, and that $d\geq s(\log n)^{60}$ (otherwise $G$ itself is a subgraph of the claimed order with minimum degree at least $s$). 
    
    Take a random subset $V'$ of $V(G)$ by keeping each vertex independently with probability $p=(\log n)^{-30}$. By the Chernoff bound, with probability $1-o(1)$, all vertices in $V(G)$ have at least $d(\log n)^{-60}$ neighbours in $V'$, and we have $n/(\log n)^{40}\leq |V'|\leq n$. Thus, if this happens, $G'=G[V']$ is a subgraph with at least $n(\log n)^{-40}\geq n^{1/3}$ vertices and minimum degree at least $d(\log n)^{-60}$. In particular, note that the minimum degree is at least $s$, and we have $1/\log n\geq 1/(\log |V'|)^2$. Therefore, we can apply \cref{thm: mu proportion} to $G'$, to deduce that there is some positive integer
    \begin{align*}
        t_0&=O_{s}\left(\max\left\{\frac{n}{(d(\log n)^{-60})^{s/(s-2)}}(\log n)^{60},\log n\right\}\right)\\
        &=O_{s}\left(\max\left\{\frac{n}{d^{s/(s-2)}}(\log n)^{180},\log n\right\}\right)
    \end{align*}
    such that, with probability $1-o(1)$, $G'$ has a (non-empty) induced subgraph on at most $t_0$ vertices in which the proportion of vertices of degree less than $s$ is at most $1/\log n$.

For each $1\leq t\leq t_0$, let $M_t$ be the number of induced subgraphs of $G$ on $t$ vertices in which the proportion of vertices of degree less than $s$ is at most $1/\log n$. Then, the expected number of such $t$-vertex graphs which are contained in $G'$ is $M_t p^{t}=M_t(\log n)^{-30 t}$. It follows that 
$\sum_{t=1}^{t_0}M_t (\log n)^{-30 t}\geq 1/2$. Since we have $\sum_{t=1}^{\infty} 3^{-t}=1/2$, it follows that there is some $1\leq t\leq t_0$ with $M_t (\log n)^{-30 t}\geq 3^{-t}$ and hence $M_t\geq (\log n)^{15t}$. 
%and hence there is some $t\leq t_0$ with $M_t\geq (\log n)^{(k/2) t}/(2n)$. Note that the statement of our lemma follows if $t<1/\mu$ -- indeed, in this case, $\mu t<1$ and so every $t$-vertex graph with at most $\mu t$ vertices of degree less than $s$ has in fact no vertex of degree less than $s$. On the other hand, if $t\geq 1/\mu \geq \log n/\log \log n$, then we have $2n\leq  (\log n)^{(k/4)t}$, giving $M_t\geq (\log n)^{(k/4)t}$.	
\end{proof}

Next, using a very similar argument, we prove the counting version of Theorem \ref{thm:almost all degree s unpruned}.

\begin{lemma} \label{lem:many with almost all degree s odd}
    Let $s\geq 5$ be an odd integer. Then for any $\eps,\mu>0$, there exist some $t_0=t_0(s,\mu)$ and $n_0=n_0(s,\eps,\mu)$ such that for each $n\geq n_0$, every $n$-vertex graph with average degree at least $n^{1-2/s+\eps}$ contains at least $n^{\eps t/2}$ induced subgraphs on $t$ vertices in which the proportion of vertices of degree less than $s$ is at most~$\mu$, for some $1\leq t\leq t_0$.
\end{lemma}

\begin{proof}
    By Theorem \ref{thm:almost all degree s unpruned} (applied with $\eps/4$ in place of $\eps$), there exist some $t_0=t_0(s,\mu)$ and $n_0'=n_0'(s,\eps,\mu)$ such that for each $n\geq n_0'$, every $n$-vertex graph with average degree at least $n^{1-2/s+\eps/4}$ contains a (non-empty) induced subgraph on at most $t_0$ vertices in which the proportion of vertices of degree less than $s$ is at most $\mu$. Let $n_0\geq (n_0')^2$ and let $n_0$ be sufficiently large compared to $s$ and $\eps$. We now show that these choices are suitable.

    Indeed, let $n\geq n_0$ and let $G$ be an $n$-vertex graph with average degree at least $n^{1-2/s+\eps}$. Take a random subset $V'$ of $V(G)$ by keeping each vertex independently with probability $p=n^{-2\eps/3}$. By standard concentration results, with probability $1-o(1)$, $G[V']$ has average degree at least $n^{1-2/s+\eps/4}$ and we have $n_0'\leq n^{1/2}\leq |V'|\leq n$. If this happens, then by the first paragraph, $G':=G[V']$ contains a (non-empty) induced subgraph on at most $t_0$ vertices in which the proportion of vertices of degree less than $s$ is at most $\mu$. We deduce that with probability $1-o(1)$, $G'$ contains such a subgraph.

    For each $1\leq t\leq t_0$, let $M_t$ be the number of induced subgraphs of $G$ on $t$ vertices in which the proportion of vertices of degree less than $s$ is at most $\mu$. Then, the expected number of such $t$-vertex graphs which are contained in $G'$ is $M_t p^{t}=M_tn^{-(2\eps/3) t}$. It follows that 
    $\sum_{t=1}^{t_0}M_t n^{-(2\eps/3) t}\geq 1/2$. Since we have $\sum_{t=1}^{\infty} 3^{-t}=1/2$, it follows that there is some $1\leq t\leq t_0$ with $M_t n^{-(2\eps/3) t}\geq 3^{-t}$ and hence $M_t\geq n^{(\eps/2)t}$.
\end{proof}

\subsection{Putting the copies together}

In this section, we build one small subgraph with large minimum degree from our large collection of small subgraphs with few low-degree vertices. 

We need some preparation for this proof. Given a collection of sets $F_1, \dots, F_r$, we say that they form a \emph{sunflower} (with $r$ \emph{petals}) if any two of the sets have the same intersection, i.e. $F_i\cap F_j=\bigcap_{\ell=1}^r F_\ell$ for all $1\leq i<j\leq r$. The set $\bigcap_{\ell=1}^r F_\ell$ is called the \emph{kernel} of the sunflower.
We will make use of the following version of the Sunflower Lemma.
\begin{theorem}[Sunflower Lemma~\cite{alweiss2021improved,Rao2020Sunflowers}]\label{thm:sunflower}
    There exists an absolute constant $C$ such that the following statement holds. Let $r\geq 2$ and $k$ be positive integers, let $X$ be any set, and let $\mathcal{F}$ be a family of subsets of $X$ such that every $F\in \mathcal{F}$ has size at most $k$. If $|\mathcal{F}|\geq (Cr\log(rk))^k$, then $\mathcal{F}$ contains $r$ sets which form a sunflower.
\end{theorem}

We will produce a small subgraph of large minimum degree via the following lemma.

\begin{lemma}\label{lem:many almost implies proper copy}
	Let $s\geq 3$ be a positive integer, let $G$ be an $n$-vertex graph, and assume that $\mu>0$ and $t$ is a positive integer such that $G$ contains $M_t\geq n^{2\mu t}(\log n)^{3t}$ induced subgraphs on $t$ vertices in which the proportion of vertices of degree less than $s$ is at most $\mu$. Then, if $n$ is large enough in terms of $s$, $G$ has a subgraph on at most $(s+1)t$ vertices in which all vertices have degree at least $s$.
\end{lemma}

\begin{proof}
By using our assumption on $t$-vertex subsets, there is some $U\subseteq V(G)$ of size at most $\mu t$ such that the vertices of degree lower than $s$ are exactly the vertices in $U$ for at least $\frac{M_t}{\binom{n}{\leq \mu t}}\geq (\log n)^{3t}$ induced subgraphs of order $t$; let $\HH$ denote the family of such subgraphs.
Our next goal is now to show that $G$ contains $s+1$ distinct induced subgraphs $H_1,\dots,H_{s+1}$ which satisfy the following properties.
\begin{enumerate}
	\item The sets $V(H_i)$ form a sunflower;
	\item If a vertex $v$ is contained in the kernel $\bigcap_i V(H_i)$, then it either has the same degree in each $H_i$ or different degrees in all the $H_i$;
	\item If a vertex $v$ has degree less than $s$ in some $H_i$ (with $v\in V(H_i)$), then $v$ is in the kernel $\bigcap_j V(H_j)$.
\end{enumerate}
To prove this claim, we form the following auxiliary set system $\mathcal{F}$. The ground set of $\mathcal{F}$ is $X=V(G)\cup (V(G)\times[0,n])$. For each subgraph $H\in\HH$, $\mathcal{F}$ contains a subset $F_H$ of $X$ given by $F_H=V(H)\cup \{(v,\deg_H(v)):v\in V(H)\}$. Thus, the family $\mathcal{F}$ contains at least $ (\log n)^{3t}$ sets. Each set in $\mathcal{F}$ contains $2t$ elements of the ground set. Since $|\mathcal{F}|\geq (\log n)^{3t}$ and $t\leq n$, if $n$ is large enough, then $|\mathcal{F}|\geq (C(s+1)\log (2(s+1)t))^{2t}$, where $C$ is the absolute constant from \cref{thm:sunflower}. Thus, \cref{thm:sunflower} gives that $\mathcal{F}$ contains a sunflower with $s+1$ petals, say $F_{H_1},\dots,F_{H_{s+1}}$. Then $H_1,\dots,H_{s+1}$ satisfy the three conditions above. (For the last property, we use that for each $i$ the set of vertices of degree less than $s$ in $H_i$ is exactly $U$.)

It remains to prove that if we have $H_1,\dots, H_{s+1}$ satisfying the three properties above, then the union of these subgraphs contains a non-empty subgraph with minimum degree at least $s$. We prove this by induction on the number of vertices in $\bigcup_i V(H_i)$. For the base case, if there is only one vertex, then (since the subgraphs are pairwise distinct) we clearly get a contradiction. Now assume there are at least $2$ vertices. If there is a vertex $v$ which does not have degree at least $s$ in the union, then it must be in the kernel, by the third property. Clearly, its degrees in the $H_i$ cannot be all different (otherwise in the union it would have degree at least $s$). Thus, by the second property, $v$ has the same degree in each $H_i$. Since the subgraphs $H_i$ are induced, $v$ has the same number of neighbours in the kernel in each $H_i$, so it must have the same number of neighbours outside the kernel in each $H_i$. If this number is nonzero, $v$ has degree at least $s+1$ in the union. Hence, all the edges leaving $v$ in some $H_i$ are in the kernel of the sunflower. For each $i\in [s+1]$, let $H_i'$ be obtained by removing $v$ and all its edges from $H_i$. Then, the graphs $H_i'$ satisfy the three properties above, since we removed the same vertex and same set of edges from each $H_i$, and if a vertex had one of its edges removed, then it must have belonged to the kernel. Thus, we may apply the induction hypothesis to $H_1',\dots,H_{s+1}'$, and the proof is complete.
\end{proof}

We are now ready to prove our first main result.

\begin{proof}[Proof of \cref{thm:even min degree}]
	We may assume that $n$ is large enough in terms of $s$, otherwise we can take the subgraph to be $G$ itself. Write $\mu=1/\log n$. By \cref{lem:many almost min degree s}, we get, for some $t=O_s(\max\{n(\log n)^{180}/d^{\frac{s}{s-2}},\log n\})$ with $t>0$, either a subgraph on $t$ vertices with minimum degree at least $s$, or $M_t\geq (\log n)^{15t}$ induced subgraphs on $t$ vertices in which the proportion of vertices of degree less than $s$ is at most $\mu$. In the former case, we are immediately done, whereas in the latter case we have $M_t\geq n^{2\mu t}(\log n)^{3t}$, so the result follows using \cref{lem:many almost implies proper copy}.
\end{proof}

Next, we prove Theorem \ref{thm:small odd subgraph}.

\begin{proof}[Proof of Theorem \ref{thm:small odd subgraph}]
    Let $\mu=\eps/5$. By Lemma \ref{lem:many with almost all degree s odd}, there exist some $t_0=t_0(s,\eps)$ and $n_0=n_0(s,\eps)$ such that for each $n\geq n_0$, every $n$-vertex graph with average degree at least $n^{1-2/s+\eps}$ contains at least $n^{\eps t_1/2}$ induced subgraphs on $t_1$ vertices in which the proportion of vertices of degree less than $s$ is at most $\mu$, for some $1\leq t_1\leq t_0$. Note that, by increasing $n_0$ if necessary, we have $n^{\eps t_1/2}\geq n^{2\mu t_1}(\log n)^{3t_1}$. Hence, by Lemma \ref{lem:many almost implies proper copy}, there is a subgraph on at most $(s+1)t_1\leq (s+1)t_0$ vertices in which all vertices have degree at least $s$. Thus, we may choose $t=(s+1)t_0$ and the theorem holds.
\end{proof}

\section{Even $s$} \label{sec:even}

In this section we prove \cref{thm: mu proportion}, which is recalled below.
\restatemuproportion*

Let us first give a proof sketch with some more details than in Section 1. For the purposes of this sketch, let us restrict to the range $d=O(n^{(s-2)/s})$. As mentioned in Section 1, we would ideally like to embed a graph obtained by first taking an $(s/2)$-blowup $H_0$ of a long path, and then replacing each vertex of $H_0$ by an independent set of size $t$, and each edge of $H_0$ by an arbitrary perfect matching between the corresponding independent sets, see Figure~\ref{fig: even case}. To do this, we will (for an appropriately chosen value of $t$) define an auxiliary graph $\AA$ on $(s/2)$-tuples of (not necessarily disjoint) $t$-element sets, where we will have an edge between two such $(s/2)$-tuples precisely when there is a perfect matching between each ``crossing'' pair of $t$-element sets, see Figure~\ref{fig:AA}. (It will also be convenient to think of this as a weighted graph, but each edge will have reasonably small weight, so for the purposes of this discussion let us ignore the weights). Then, the ideal graph described above essentially corresponds to a long path $\mathcal{P}$ in this auxiliary graph $\AA$. However, this is not entirely true, as there may be vertices in $G$ appearing multiple times in $\mathcal{P}$ due to overlaps (either within an $(s/2)$-tuple, or across different $(s/2)$-tuples).

To solve the issue of overlapping vertices (and edges), first observe that they only really cause problems for us if some overlaps appear among the neighbours of a vertex (as otherwise we still have a total degree of at least $s$ for each vertex in the interior of $\mathcal{P}$). We will obtain $\mathcal{P}$ by taking a random walk on the auxiliary graph $\AA$, starting from the stationary distribution. Then, our goal is to show that each vertex will have, with high probability, at least $s$ neighbours among the vertices on this path (i.e. few overlaps among its neighbours). To do so, we show that if a vertex has, with reasonably high probability, two of its $s$ matched neighbours coinciding, then, also with reasonably high probability, \emph{all} of its $s$ matched vertices will coincide. This is true because, conditioning on visiting $\XX$ in $\mathcal{P}$ (so that $\XX$ is an $(s/2)$-tuple of $t$-sets), the neighbours of a vertex inside $\XX$ can be shown to be independent and identically distributed, using the properties of the stationary distribution. Therefore, it is enough to show that it is unlikely that \emph{many} vertices have the property that \emph{all} of their neighbours coincide. This will be verified using a double counting argument: many such overlaps would collapse the collection of perfect matchings into a smaller structure, and we can show that such smaller structures cannot have too many copies appearing.\medskip

Let us now proceed with the formal proof. Note that the following result of O.~Janzer, implicit in~\cite[Theorem 1.14]{janzer2023rainbow}, immediately implies the statement of \cref{thm:even min degree} in the range when $d=\Omega(n^{(s-2)/s}(\log n)^4)$. 

\begin{theorem}[O. Janzer \cite{janzer2023rainbow}]\label{thm:dense case even}
    For every even integer $s\geq 4$ there exists a constant $C_s$ such that whenever $G$ is an $n$-vertex graph with at least $C_s n^{2-2/s}(\log n)^{14/s}$ edges, then $G$ has an $s$-regular subgraph on $O_s(\log n)$ vertices. 
\end{theorem}
We note that when $d=\Omega(n^{(s-2)/s+\delta})$ for some $\delta>0$, then the $s$-regular subgraph can even be guaranteed to have order $O_{s,\delta}(1)$ \cite[Theorem 1.15]{janzer2023rainbow}.

The following standard lemma will be used several times. (It can be proved by repeatedly removing a vertex of minimum degree from our graph until the conclusion is satisfied.)

\begin{lemma}\label{lem:min+averagedegree}
    Let $G$ be a graph with minimum degree $d$ and average degree $\bar{d}$. Then $G$ has an induced subgraph whose minimum degree $d'$ is at least $d$, its average degree $\bar{d'}$ is at least $\bar{d}$, and $\bar{d'}\leq 2d'$.
\end{lemma}

As explained above, the proof of Theorem \ref{thm: mu proportion} will involve counting certain structures consisting of perfect matchings, for which the following two lemmas will be useful.

 \begin{lemma}\label{lem:matching count}
        Assume that $s\geq 4$ is even, $G$ is a graph, and $t>0$ such that $G$ does not have a non-empty subgraph of order at most $2t$ with minimum degree at least $s$. Assume that $\alpha, \beta\in (0,1]$, and let $X, Y$ be disjoint subsets of $V(G)$ of orders $\alpha t, \beta t$, respectively. Then, the number of matchings between $X$ and $Y$ saturating $Y$ is at most $(\frac{\alpha+\beta}{\beta}s)^{\beta t}$. In particular, if $|X|=|Y|=t$, then the number of perfect matchings between $X$ and $Y$ is at most $(2s)^{t}$.
    \end{lemma}
    \begin{proof}
    	Note that, by \cref{lem:min+averagedegree}, the average degree in $G[X\cup Y]$ must be less than $2s$. Thus, the number of matchings between $X$ and $Y$ saturating $Y$ is at most $\prod_{y\in Y}\deg_{X}(y)\leq (\frac{1}{\beta t}\sum_{y\in Y} \deg_{X}(y))^{\beta t}\leq (\frac{\alpha+\beta}{\beta}s)^{\beta t}$, where we used the inequality between the arithmetic and geometric means in the first step.
    \end{proof}

\begin{lemma}\label{lem:count of matching "blowups"}
        Assume that $s\geq 4$ is even, $G$ is a graph of $n$ vertices and minimum degree at least $d\geq 100s$, and $100s\leq t\leq n/(100s)$ such that $G$ does not have a subgraph of order at most $20st$ with minimum degree at least $s$. Let $W\subseteq V(G)$ be a set of size at most $8st$, assume that $n$ is large enough, and let $\alpha\in (0,1]$ such that $\alpha t$ is an integer. Then, the number of $s$-tuples of sets $S_1, \dots, S_{s/2}, T_1,\dots, T_{s/2}$ of size $\alpha t$ each such that for each $i, j\in [s/2]$, $S_i$, $T_j$ and $W$ are pairwise disjoint, and there is a perfect matching between $S_i$ and $T_j$, is at least
        $$\left(\Omega_{s}(1)\cdot  \frac{d^{s}}{(n/(\alpha t))^{s-4}}\right)^{\alpha ts/4}.$$

    \end{lemma}
	\begin{proof}
	First note that for all sets $U\subseteq V(G)$ of size $|U|\leq 10st$, the number of vertices which have at least $4s$ neighbours in $U$ is at most $10st$, by \cref{lem:min+averagedegree}. 

    Define an auxiliary graph $\MM$ by letting its vertex set consist of the $(\alpha t)$-sized sets of vertices in $V(G)\setminus W$, and two such $(\alpha t)$-sets are joined by an edge precisely when they are disjoint and there is a perfect matching between them. Notice that the quantity we wish to lower bound is exactly the number of graph homomorphisms from $K_{s/2, s/2}$ to $\MM$. 
    
	Let us count the number of matchings of size $\alpha t$ in $G[V(G)\setminus W]$. Select the edges appearing in the matching one-by-one, always choosing a vertex which has less than $4s$ neighbours in $W$  and the already used vertices combined, and choose one of its at least $d-4s$ ``unused'' neighbours to be matched to it. Using the observation in the first paragraph, and noting that each matching is counted at most $2^{\alpha t}(\alpha t)!$ times, we see that there are at least $\frac{1}{2^{\alpha t} (\alpha t)!}((n-20st)(d-4s))^{\alpha t}\geq \left(\frac{nd}{4\alpha t}\right)^{\alpha t}$ such matchings of size $\alpha t$.

    Using Lemma~\ref{lem:matching count}, it follows that the number of edges of the auxiliary graph $\MM$ is at least $$|E(\MM)|\geq \frac{1}{2}\left(\frac{nd}{8s\alpha t}\right)^{\alpha t}.$$

    It is well known that Sidorenko's conjecture holds for complete bipartite graphs. Therefore, the number of homomorphic copies of $K_{s/2,s/2}$ in $\MM$ is at least
    \begin{align*}
        |V(\MM)|^{s}\left(\frac{2|E(\MM)|}{|V(\MM)|^2}\right)^{(s/2)^2}&=\left(|V(\MM)|^{2-s}\left(2|E(\MM)|\right)^{s/2}\right)^{s/2}\\
        &\geq \left(\frac{1}{\binom{n}{\alpha t}^{s-2}}\left(\frac{nd}{8s\alpha t}\right)^{\alpha ts/2}\right)^{s/2}\\
        &\geq \left(\left(\frac{en}{\alpha t}\right)^{-\alpha t(s-2)}\left(\frac{nd}{8s\alpha t}\right)^{\alpha ts/2}\right)^{s/2}\\
        &=\left(\frac{1}{e^{2(s-2)}(8s)^s}\cdot d^s\cdot (n/(\alpha t))^{4-s}\right)^{\alpha ts/4},
    \end{align*}
giving the claimed bound.
    \end{proof}
    
We are ready to prove the main result of this section.

\begin{proof}[Proof of \cref{thm: mu proportion}]
	We may assume that $n$ is sufficiently large (in terms of $s$), and also that $\frac{n}{d^{s/(s-2)}}(\log n)^{60}\leq n$, otherwise we can take the subgraph to be $G$ itself. Moreover, as noted before, by \cref{thm:dense case even}, we may also assume $d\leq n^{\frac{s-2}{s}}(\log n)^{4}$ and hence $\frac{n}{d^{s/(s-2)}}(\log n)^{60} \geq (\log n)^{50}$. In particular, writing
	$$t=\left\lfloor \frac{n}{d^{s/(s-2)}} (\log n)^{45} \right\rfloor,$$
	we have
	$$(\log n)^{30}\leq t\leq n/(\log n)^{15},$$
	and the claimed bound on the order of the subgraph we can find is $O_{s}(t (\log n)^{15})$. From now on, we may assume that $G$ does not have a subgraph of order $O_s(t)$ with minimum degree at least $s$. Let us also write
    $$\mu=\frac{1}{(\log n)^2}$$
    for the proportion of vertices that are allowed to have degree lower than $s$.
    
	Define an auxiliary graph $\MM$ for which the vertex set is $\binom{V(G)}{t}$, and $X \in V(\MM)$ is joined to $Y \in V(\MM)$ precisely when $X$ and $Y$ are disjoint and there is a perfect matching between $X$ and $Y$. Moreover, we view this as a weighted graph, where an edge between $X$ and $Y$ has weight equal to  the number of perfect matchings between $X$ and $Y$. Note that, by \cref{lem:matching count}, each edge of $\MM$ has weight at most $(2s)^t$.

    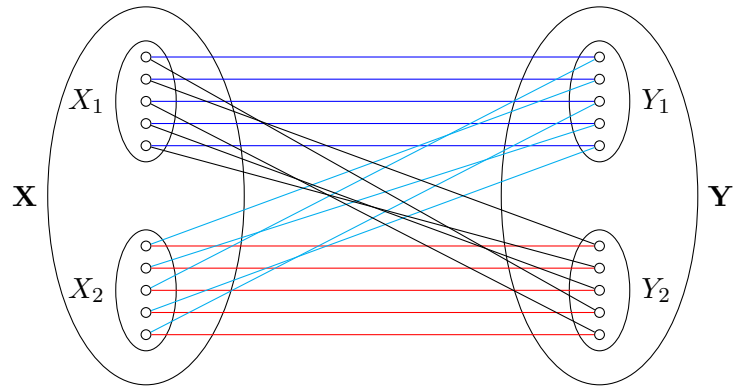
\begin{figure}[b]
    \centering

\begin{tikzpicture}[thin,
	every node/.style={draw,circle},
    every fit/.style={ellipse,draw,text width=1cm}	]

	% the vertices of X_{1}
	\begin{scope}[xshift=0cm,yshift=0cm,start chain=going below,node distance=1.5mm, inner sep=1.3pt]
		\foreach \i in {0,1,2,3,4}
		\node[%fsnode,
		on chain] (x1\i) {};
	\end{scope}
	% the vertices of X_{2}
	\begin{scope}[xshift=0cm,yshift=2.5cm,start chain=going below,node distance=1.5mm, inner sep=1.3pt]
		\foreach \i in {0,1,2,3,4}
		\node[%fsnode,
		on chain] (x2\i) {};
	\end{scope}
    
    % the vertices of Y_{1}
	\begin{scope}[xshift=6cm,yshift=0cm,start chain=going below,node distance=1.5mm, inner sep=1.3pt]
		\foreach \i in {0,1,2,3,4}
		\node[%fsnode,
		on chain] (y1\i) {};
	\end{scope}
	% the vertices of Y_{2}
	\begin{scope}[xshift=6cm,yshift=2.5cm,start chain=going below,node distance=1.5mm, inner sep=1.3pt]
		\foreach \i in {0,1,2,3,4}
		\node[%fsnode,
		on chain] (y2\i) {};
	\end{scope}
	
	%labels shown here are not consistent with tikz names
\coordinate (c) at ($(x10)!0.5!(x14)$);
\draw (c) ellipse [x radius=4mm, y radius=8mm];
\node[left=3.1mm, draw=none] at (c) {$X_2$};
\coordinate (c) at ($(x20)!0.5!(x24)$);
\draw (c) ellipse [x radius=4mm, y radius=8mm];
\node[left=3.1mm, draw=none] at (c) {$X_1$};
\coordinate (c) at ($(x10)!0.5!(x24)$);
\draw (c) ellipse [x radius=13mm, y radius=25mm];
\node[left=12mm, draw=none] at (c) {$\XX$};

\coordinate (c) at ($(y10)!0.5!(y14)$);
\draw (c) ellipse [x radius=4mm, y radius=8mm];
\node[right=3.1mm, draw=none] at (c) {$Y_2$};
\coordinate (c) at ($(y20)!0.5!(y24)$);
\draw (c) ellipse [x radius=4mm, y radius=8mm];
\node[right=3.1mm, draw=none] at (c) {$Y_1$};
\coordinate (c) at ($(y10)!0.5!(y24)$);
\draw (c) ellipse [x radius=13mm, y radius=25mm];
\node[right=12mm, draw=none] at (c) {$\YY$};

\def\colors{{"red","blue","cyan","black","cyan","brown","gray","black"}}

    \foreach \i in {1,2}
    \foreach \k in {0,1,2,3,4}
    \pgfmathparse{\colors[\i-1]}
    \let\edgecolor\pgfmathresult
    \draw[\edgecolor] (x\i\k) to (y\i\k);

    \pgfmathparse{\colors[2]}
    \let\edgecolor\pgfmathresult
    \draw[\edgecolor] (x10) to (y21);
    \draw[\edgecolor] (x11) to (y23);
    \draw[\edgecolor] (x12) to (y20);
    \draw[\edgecolor] (x13) to (y24);
    \draw[\edgecolor] (x14) to (y22);

    \pgfmathparse{\colors[3]}
    \let\edgecolor\pgfmathresult
    \draw[\edgecolor] (x20) to (y13);
    \draw[\edgecolor] (x21) to (y10);
    \draw[\edgecolor] (x22) to (y14);
    \draw[\edgecolor] (x23) to (y12);
    \draw[\edgecolor] (x24) to (y11);

\end{tikzpicture}
\captionsetup{justification=centering}
\caption{An illustration of an edge $\XX\YY$ of the auxiliary graph $\mathcal{A}$ in the case $s=4$, $t=5$. The edge between $\XX=(X_1,X_2)$ and $\YY=(Y_1,Y_2)$ is obtained by a perfect matching between each of the pairs $X_i, Y_j$.}    
\label{fig:AA}
\end{figure}

	Define another auxiliary graph $\AA$ for which the vertex set is the collection of $(s/2)$-tuples $\XX=(X_1,\dots,X_{s/2})$ where $X_1,\dots,X_{s/2}$ are (not necessarily disjoint) subsets of $V(G)$ of size $t$ each, and $(X_1,\dots,X_{s/2})$ is joined to $(Y_1,\dots,Y_{s/2})$ precisely when $X_iY_j$ is an edge of $\MM$ for each $i,j$ (see Figure~\ref{fig:AA}). Moreover, we view this as a weighted graph, where an edge between $(X_1,\dots,X_{s/2})$ and $(Y_1,\dots,Y_{s/2})$ has weight equal to the product of $\MM$-weights of the edges $X_iY_j$. By the above, each edge in $\AA$ has weight at most $(2s)^{(s/2)^2t}\leq (2s)^{s^2 t}$. Given $\XX$ as above, by the vertices in $\XX$ we will mean the vertices in $\bigcup_{i=1}^{s/2} X_i$. For all large enough $n$, we get that $e(\AA)\not=0$ by \cref{lem:count of matching "blowups"} with $\alpha=1$ and $W=\emptyset$.

	Let $$N=\lceil (3s+1)/\mu ^2\rceil,$$ and take a random walk on $\AA$ starting from the stationary distribution, say $\XX_{-N},\dots,\XX_N$. That is, for each $\XX\in V(\AA)$, $\PP(\XX_{-N}=\XX)=\frac{\deg_\AA(\XX)}{2e(\AA)}$ (with $\deg_\AA(\XX)$ and $e(\AA)$ defined as the corresponding sums of weights), and given $\XX_{-N},\dots,\XX_i$, the next step $\XX_{i+1}$ is a random neighbour of $\XX_i$, with probabilities proportional to the edge weights from $\XX_{i}$. Note that this implies that each $\XX_i$ is also distributed according to the stationary distribution.  Moreover, the random walk is reversible, i.e. given any $\XX_i,\XX_{i+1},\dots,\XX_N$ for some $i>-N$, the preceding vertex $\XX_{i-1}$ is also a random $\AA$-neighbour of $\XX_i$, with probabilities proportional to the corresponding edge weights.  We will show that (if $n$ is large enough), with positive probability, we are able to find a subset of the set of the vertices appearing in this walk that span a subgraph in which the proportion of vertices of degree less than $s$ is at most $\mu$. Once we establish this, the proof will be complete, as the number of vertices appearing in the random walk is at most $(2N+1)\cdot (s/2)\cdot t=O_{s}(t/\mu ^2)=O_s(t(\log n)^{15})$.\medskip
	
	Let $\eta>0$ be such that $$\eta\leq \frac{\mu}{6Ns^3},\textnormal{\hspace{12pt}}\eta=\Theta_s(\mu^3)\textnormal{\hspace{12pt} and \hspace{12pt}}\eta t\in\mathbb{Z}.$$ Note that $\mu^3= (\log n)^{-6}$ and $t\geq (\log n)^{30}$, so it is possible to pick such an $\eta$, and we have
	$$\eta t=\Omega_{s}((\log n)^{15}).$$
		Given $\XX=(X_1,\dots,X_{s/2})\in V(\AA)$ with $\deg_\AA(\XX)>0$, we say that $\XX$ is \emph{bad} if, picking $\YY^{(1)},\YY^{(2)}$ independently at random from the $\AA$-neighbours of $\XX$, with the probabilities $\PP(\YY^{(i)}=\YY)$ proportional to the edge weights $\XX \YY$, then the probability that	at least $s^3\eta t$ vertices in $\XX$ have fewer than $s$ neighbours in the union of vertices in $\YY^{(1)}$ and $\YY^{(2)}$
is at least $\frac{1}{6N}$. We note that this notion will be useful because in the random walk above, $\XX_{i-1}$ and $\XX_{i+1}$ are independent random neighbours of $\XX_i$ chosen according to the distribution above, so if $\XX$ is not bad, then we expect most of its vertices to have degree at least $s$ in the vertices appearing in $\XX_{i-1},\XX_{i+1}$. The following claim says that if $\XX$ is bad, then a ``large'' proportion of the $\AA$-neighbours $\YY$ of $\XX$ contain some specific set $Y^\star$ of size $\eta t$ in each $Y_i$ (causing the overlaps which make the degrees to be lower than $s$). We will then use this claim to deduce that most $\XX_i$ appearing in the random walk are good.

\begin{claim}\label{claim:bad X}
    Assume that $\XX=(X_1,\dots,X_{s/2})\in V(\AA)$ is bad. Then there exist $X_i^\star\subseteq X_i$ for all $i\in[s/2]$, as well as  $Y^\star\subseteq V(G)$, each of size $\eta t$, together with perfect matchings $N_i$ between $X_i^\star$ and $Y^\star$, such that there are $(\Omega_{s}(\eta))^{\eta ts} \deg_{\AA}(\XX)^{\frac{1}{s/2}}$ tuples $(Y,(M_i)_{i=1}^{s/2})$ such that
    \begin{itemize}
    	\item $Y$ is a subset of $V(G)$ of size $t$ with $Y^\star\subseteq Y$, and
    	\item for all $i\in[s/2]$, $M_i$ is a perfect matching between $X_i$ and $Y$ such that $N_i\subseteq M_i$.
    \end{itemize}
\end{claim}
\begin{proof}
	Let $\mathcal{S}$ be the set of tuples $\left(Z,(M_i)_{i=1}^{s/2}\right)$ such that $Z$ is a subset of $V(G)\setminus \bigcup_i X_i$ of size $t$, and for each $i$, $M_i$ is a perfect matching between $X_i$ and $Z$. Note that the number of times a set $Z$ appears (as a first coordinate) in $\mathcal{S}$ is equal to the product of the $\MM$-weights of $X_iZ$. Hence, if $\left(Z_1,(M_{i,1})_{i=1}^{s/2}\right), \dots, \left(Z_{s/2},(M_{i,s/2})_{i=1}^{s/2}\right)$ are $s/2$ independent, uniformly random elements of $\mathcal{S}$, then $(Z_1,\dots,Z_{s/2})$ is a random $\AA$-neighbour of $\XX$ chosen according to the $\AA$-edge weights. Moreover, $|\mathcal{S}|=(\deg_\AA(\XX))^{\frac{1}{s/2}}$.

Thus, according to the definitions above, $\XX$ is bad if, picking $s$ elements $\left(Z_1,(M_{i,1})_{i=1}^{s/2}\right), \dots, \left(Z_{s},(M_{i,s})_{i=1}^{s/2}\right)$ of $\mathcal{S}$ uniformly at random, we have that the probability that at least $s^3\eta t$ vertices in $\XX$ have less than $s$ neighbours in $Z_1\cup\dots\cup Z_s$
is at least $\frac{1}{6N}$. Note that, if the number of such vertices with fewer than $s$ neighbours is at least $s^3\eta t$, then there is some $i^\star\in[s/2]$ and $X'\subseteq X_{i^\star}$ of size $s^2\eta t$ such that each vertex in $X'$ has less than $s$ neighbours in $Z_1\cup\dots\cup Z_s$. However, if a vertex in $X_i$ has degree less than $s$ in $Z_1\cup\dots\cup Z_s$, then it must have the same neighbour in two of its $s$ matchings $M_{i,1},\dots,M_{i,s}$. Hence, there is some $X^\star\subseteq X'$ of size $\eta t$ and $j_1^\star,j_2^\star\in[s]$ distinct such that every vertex in $X^\star$ has the same neighbour in the two matchings $M_{i^\star, j_1^\star}$ and $M_{i^\star, j_2^\star}$.

 Moreover, if $\XX$ is bad, then at least $\frac{|\mathcal{S}|^s}{6N(s/2)\binom{t}{\eta t}s^2}$ $s$-tuples from $\mathcal{S}$ give the same $i^\star, X^\star, j_1^\star, j_2^\star$. Indeed, as $\XX$ is bad, the proportion of the $|\mathcal{S}|^s$ tuples having a corresponding $(i^\star, X^\star,j_1^\star,j_2^\star)$ is at least $1/(6N)$, and $i^\star$ can take $s/2$ different values, after which $X^\star$ can take $\binom{t}{\eta t}$ values and $j_1^\star, j_2^\star$ can take at most $s$ values each. Fix such an $(i^\star, X^\star,j_1^\star,j_2^\star)$ corresponding to at least $\frac{|\mathcal{S}|^s}{3Ns^3\binom{t}{\eta t}}$ $s$-tuples, and for each (injective) function $f: X^\star \to V(G)$, let $a_f$ be the number of $\left(Z,(M_i)_{i=1}^{s/2}\right)\in \mathcal{S}$ such that $f(X^\star)\subseteq Z$ and for all $x\in X^\star$, $x$ is matched to $f(x)$ in $M_{i^\star}$. Then, counting the number of ways we can have each element of $X^\star$ matched to the same element in the two matchings $M_{i^\star,j_1^\star}$ and $M_{i^\star,j_2^\star}$, we have $\sum_f a_f^2|\mathcal{S}|^{s-2}\geq\frac{|\mathcal{S}|^s}{3Ns^3\binom{t}{\eta t}}$. However, $a_f^2\leq a_f\max_{g} a_{g}$ and $\sum_f a_f=|\mathcal{S}|$, giving $\max_{g}a_{g}\geq \frac{|\mathcal{S}|}{3Ns^3\binom{t}{\eta t}}$. Let $f^\star$ be a choice of $f$ such that $a_{f^\star}$ is maximal, so $a_{f^\star}\geq  \frac{|\mathcal{S}|}{3Ns^3\binom{t}{\eta t}}$, and let $Y^\star= f^\star(X^\star)$. Note that, by \cref{lem:matching count}, for each $i\in[s/2]$, the number of matchings between $X_i$ and $Y^\star$ saturating $Y^\star$ is at most $\left(\frac{2}{\eta}s\right)^{\eta t}$. Hence, for some matchings $N_i$ between $X_i$ and $Y^\star$ saturating $Y^\star$, the number of $\left(Z,(M_i)_{i=1}^{s/2}\right)\in \mathcal{S}$ such that $Y^\star\subseteq Z$ and for all $i$ $N_i\subseteq M_i$, is at least $\frac{|\mathcal{S}|}{3Ns^3\binom{t}{\eta t}\left(\frac{2}{\eta}s\right)^{\eta ts/2}}=\deg_\AA(\XX)^{\frac{1}{s/2}}\cdot (\Omega_s(\eta^{s/2+1}))^{\eta t}\cdot \Omega_s(1/N)=\deg_\AA(\XX)^{\frac{1}{s/2}}\cdot (\Omega_s(\eta^{s}))^{\eta t}$. Letting $X_i^\star$ be $V(N_i)\cap X_i$ for all $i$, the claim follows.
    \renewcommand{\qedsymbol}{$\blacksquare$}
\end{proof}
    \renewcommand{\qedsymbol}{$\square$}

    Using the claim above, we can now show that few edges of $\AA$ come from some bad $\XX$.
    \begin{claim}\label{claim:few bad}
        We have
        \begin{equation*}%\label{eq:few bad}
\frac{\sum_{\XX\textnormal{ bad}}\deg_\AA(\XX)}{e(\AA)}\leq e^{-(\Omega_s(\log n))^{15}}.
\end{equation*}
    \end{claim}
    \begin{proof}

    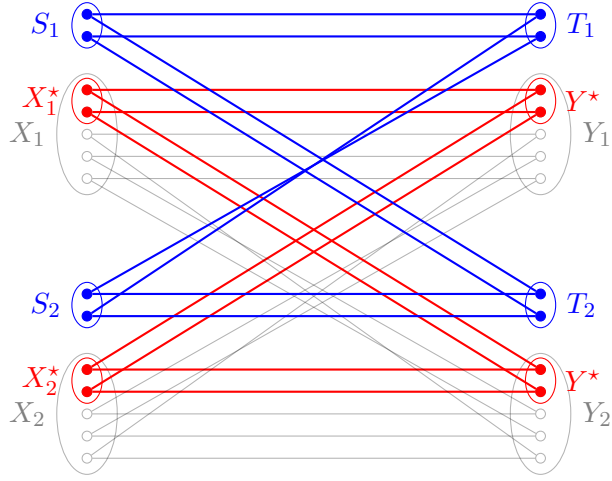
\begin{figure}[t]
    \centering
\begin{tikzpicture}[thin,
	every node/.style={draw,circle},
    every fit/.style={ellipse,draw,text width=1cm}	]

	% the vertices of X_{1}
	\begin{scope}[xshift=0cm,yshift=0cm,start chain=going below,node distance=1.5mm, inner sep=1.3pt]
		\foreach \i in {0,1,2,3,4}
		\node[%fsnode,
		on chain, opacity=0.3] (x1\i) {};
	\end{scope}
	% the vertices of X_{2}
	\begin{scope}[xshift=0cm,yshift=3.7cm,start chain=going below,node distance=1.5mm, inner sep=1.3pt]
		\foreach \i in {0,1,2,3,4}
		\node[%fsnode,
		on chain, opacity=0.3] (x2\i) {};
	\end{scope}

	% the vertices of X_{1}^\star
	\begin{scope}[xshift=0cm,yshift=0cm,start chain=going below,node distance=1.5mm, inner sep=1.3pt]
		\foreach \i in {0,1}
		\node[on chain, color=red, fill] (x1\i) {};
	\end{scope}

	% the vertices of X_{2}^\star
	\begin{scope}[xshift=0cm,yshift=3.7cm,start chain=going below,node distance=1.5mm, inner sep=1.3pt]
		\foreach \i in {0,1}
		\node[on chain, color=red, fill] (x2\i) {};
	\end{scope}

	% the vertices of S_{1}
	\begin{scope}[xshift=0cm,yshift=1cm,start chain=going below,node distance=1.5mm, inner sep=1.3pt]
		\foreach \i in {0,1}
		\node[on chain, color=blue, fill] (s1\i) {};
	\end{scope}
	% the vertices of S_{2}
	\begin{scope}[xshift=0cm,yshift=4.7cm,start chain=going below,node distance=1.5mm, inner sep=1.3pt]
		\foreach \i in {0,1}
		\node[on chain, color=blue, fill] (s2\i) {};
	\end{scope}
	% the vertices of T_{1}
	\begin{scope}[xshift=6cm,yshift=1cm,start chain=going below,node distance=1.5mm, inner sep=1.3pt]
		\foreach \i in {0,1}
		\node[on chain, color=blue, fill] (t1\i) {};
	\end{scope}
	% the vertices of T_{1}
	\begin{scope}[xshift=6cm,yshift=4.7cm,start chain=going below,node distance=1.5mm, inner sep=1.3pt]
		\foreach \i in {0,1}
		\node[on chain, color=blue, fill] (t2\i) {};
	\end{scope}

    % the vertices of Y_{1}
	\begin{scope}[xshift=6cm,yshift=0cm,start chain=going below,node distance=1.5mm, inner sep=1.3pt]
		\foreach \i in {0,1,2,3,4}
		\node[%fsnode,
		on chain, opacity=0.3] (y1\i) {};
	\end{scope}
	% the vertices of Y_{2}
	\begin{scope}[xshift=6cm,yshift=3.7cm,start chain=going below,node distance=1.5mm, inner sep=1.3pt]
		\foreach \i in {0,1,2,3,4}
		\node[%fsnode,
		on chain, opacity=0.3] (y2\i) {};
	\end{scope}

	% the vertices of Y^\star in Y_2
	\begin{scope}[xshift=6cm,yshift=0cm,start chain=going below,node distance=1.5mm, inner sep=1.3pt]
		\foreach \i in {0,1}
		\node[on chain, color=red, fill] (y1\i) {};
	\end{scope}

	% the vertices of Y^\star in Y_2
	\begin{scope}[xshift=6cm,yshift=3.7cm,start chain=going below,node distance=1.5mm, inner sep=1.3pt]
		\foreach \i in {0,1}
		\node[on chain, color=red, fill] (y2\i) {};
	\end{scope}
	
	%labels shown here are not consistent with tikz names
\coordinate (c) at ($(x10)!0.5!(x14)$);
\draw[, opacity=0.3] (c) ellipse [x radius=4mm, y radius=8mm];
\node[left=3.1mm, draw=none] at (c) {\textcolor{gray}{$X_2$}};

\coordinate (c) at ($(x20)!0.5!(x24)$);
\draw[, opacity=0.3] (c) ellipse [x radius=4mm, y radius=8mm];
\node[left=3.1mm, draw=none] at (c) {\textcolor{gray}{$X_1$}};

\coordinate (c) at ($(y10)!0.5!(y14)$);
\draw[, opacity=0.3] (c) ellipse [x radius=4mm, y radius=8mm];
\node[right=3.1mm, draw=none] at (c) {\textcolor{gray}{$Y_2$}};

\coordinate (c) at ($(y20)!0.5!(y24)$);
\draw[, opacity=0.3] (c) ellipse [x radius=4mm, y radius=8mm];
\node[right=3.1mm, draw=none] at (c) {\textcolor{gray}{$Y_1$}};

\coordinate (c) at ($(x20)!0.5!(x21)$);
\draw[red] (c) ellipse [x radius=2mm, y radius=3mm];
\node[left=1mm, draw=none] at (c) {\textcolor{red}{$X_1^\star$}};

\coordinate (c) at ($(x10)!0.5!(x11)$);
\draw[red] (c) ellipse [x radius=2mm, y radius=3mm];
\node[left=1mm, draw=none] at (c) {\textcolor{red}{$X_2^\star$}};

\coordinate (c) at ($(y20)!0.5!(y21)$);
\draw[red] (c) ellipse [x radius=2mm, y radius=3mm];
\node[right=1mm, draw=none] at (c) {\textcolor{red}{$Y^\star$}};

\coordinate (c) at ($(y10)!0.5!(y11)$);
\draw[red] (c) ellipse [x radius=2mm, y radius=3mm];
\node[right=1mm, draw=none] at (c) {\textcolor{red}{$Y^\star$}};

\coordinate (c) at ($(s20)!0.5!(s21)$);
\draw[blue] (c) ellipse [x radius=2mm, y radius=3mm];
\node[left=1mm, draw=none] at (c) {\textcolor{blue}{$S_1$}};

\coordinate (c) at ($(s10)!0.5!(s11)$);
\draw[blue] (c) ellipse [x radius=2mm, y radius=3mm];
\node[left=1mm, draw=none] at (c) {\textcolor{blue}{$S_2$}};

\coordinate (c) at ($(t20)!0.5!(t21)$);
\draw[blue] (c) ellipse [x radius=2mm, y radius=3mm];
\node[right=1mm, draw=none] at (c) {\textcolor{blue}{$T_1$}};

\coordinate (c) at ($(t10)!0.5!(t11)$);
\draw[blue] (c) ellipse [x radius=2mm, y radius=3mm];
\node[right=1mm, draw=none] at (c) {\textcolor{blue}{$T_2$}};

\def\colors{{"black", "black", "black","black"}}

    \foreach \i in {1,2}
    \foreach \k in {0,1}
    \draw[red, thick] (x\i\k) to (y\i\k);

    \foreach \i in {1,2}    
    \foreach \k in{2,3,4}
    \pgfmathparse{\colors[\i-1]}
    \let\edgecolor\pgfmathresult
    \draw[\edgecolor, opacity=0.3] (x\i\k) to (y\i\k);

    \pgfmathparse{\colors[2]}
    \let\edgecolor\pgfmathresult
    \draw[red, thick] (x10) to (y20);
    \draw[red, thick] (x11) to (y21);
    \draw[\edgecolor, opacity=0.3] (x12) to (y23);
    \draw[\edgecolor, opacity=0.3] (x13) to (y24);
    \draw[\edgecolor, opacity=0.3] (x14) to (y22);

    \pgfmathparse{\colors[3]}
    \let\edgecolor\pgfmathresult
    \draw[red, thick] (x20) to (y10);
    \draw[red, thick] (x21) to (y11);
    \draw[\edgecolor, opacity=0.3] (x22) to (y14);
    \draw[\edgecolor, opacity=0.3] (x23) to (y12);
    \draw[\edgecolor, opacity=0.3] (x24) to (y13);

    \draw[blue, thick] (s10) to (t10);
    \draw[blue, thick] (s11) to (t11);
    \draw[blue, thick] (s20) to (t20);
    \draw[blue, thick] (s21) to (t21);
    \draw[blue, thick] (s10) to (t21);
    \draw[blue, thick] (s11) to (t20);
    \draw[blue, thick] (s20) to (t10);
    \draw[blue, thick] (s21) to (t11);

\end{tikzpicture}
\captionsetup{justification=centerlast}
\caption{An illustration of the tuples we are double counting in the proof of Claim~\ref{claim:few bad}. For the purposes of this figure, $s=4$, $t=5$ and $\eta t=2$. Note that we require $Y^\star\subseteq Y_j$ for all $j$, and each copy of $Y^\star$ is matched to $X_i^\star$ in $X_i$ (via $N_i$, in red) for all $i$. Then $\tilde{\XX}, \tilde{\YY}$ are obtained by removing the red part (corresponding to matchings $N_i$, i.e. the repeated set $Y^\star$ and its neighbours $X_i^\star$), and instead adding arbitrary sets of sizes $\eta t$ with matchings between them (shown in blue).}
\label{fig:quadruples}
\end{figure}

	For each bad $\XX$, we fix a choice of sets and matchings $X_i^\star, Y^\star, N_{i}$ satisfying Claim~\ref{claim:bad X}.
	We double count the tuples $(\XX,\YY,(M_{i,j})_{i,j\in[s/2]},\tilde\XX,\tilde \YY,(\tilde{M}_{i,j})_{i,j\in[s/2]})$ such that the following properties hold (see Figure~\ref{fig:quadruples}):
	\begin{itemize}
	\item $\XX\YY$ and $\tilde \XX\tilde \YY$ are edges of $\AA$, and for all $i,j\in[s/2]$, $M_{i,j}$ is a perfect matching between $X_{i}$ and $Y_j$, and $\tilde{M}_{i,j}$ is a perfect matching between $\tilde{X}_{i}$ and $\tilde{Y}_j$,
	\item $\XX$ is bad,
    \item We have $Y^\star\subseteq Y_j$ for all $j$, and for all $i,j\in[s/2]$, $N_{i}\subseteq M_{i,j}$ (where $X_i^\star, Y^\star, N_{i}$ are the sets and matchings corresponding to $\XX$, coming from Claim~\ref{claim:bad X}),
	\item For all $i\in[s/2]$, $\tilde X_i=(X_i\setminus X_i^\star)\cup S_i$ for some $S_i$ of size $\eta t$,
	\item For all $j\in[s/2]$, $\tilde Y_j=(Y_j\setminus Y^\star)\cup T_j$ for some $T_j$ of size $\eta t$,
	\item For all $i,j\in[s/2]$, there is a perfect matching between $S_i$ and $T_j$. Moreover, $\tilde{M}_{i,j}$ is obtained from $M_{i,j}$ by removing $N_{i}$ and adding some perfect matching between $S_i$ and $T_j$.
	\end{itemize}  
	
	On the one hand, by choosing $(\tilde \XX,\tilde \YY, (\tilde{M}_{i,j})_{i,j\in[s/2]})$ first, then the sets $S_i$ and $T_j$ and finally the sets $X_i ^\star$ (note that $\XX$ then determines the missing part $Y^\star$ of $\YY$ and $N_i$ of each $M_{i,j}$), we see that the number of such tuples is at most
	$$2e(\AA)\cdot \binom{t}{\eta t}^{s}\binom{n}{\eta t}^{s/2}=2e(\AA)\cdot \left(O_s\left(\frac{n}{\eta^3 t}\right)\right)^{\eta ts/2}.
	$$
	
	On the other hand, consider the number of such tuples by picking $\XX$ first (then the $X_i^\star$, $Y^\star$ and $N_i$ are uniquely determined). Then, by the conclusion of Claim~\ref{claim:bad X}, we can pick $(\YY,(M_{i,j})_{i,j\in[s/2]})$ in at least $\left(\Omega_{s}(\eta)^{\eta ts} \deg_\AA(\XX)^{2/s}\right)^{s/2}
    =\Omega_{s}(\eta)^{\eta ts^2/2}{\deg_\AA(\XX)}$ ways. Finally, we pick the sets $S_i$ and $T_j$ and a matching between each pair $S_i, T_j$ (and we will also guarantee that for all $i,j$, $S_i\cap T_j=\emptyset$, and these sets are disjoint from the vertices appearing in $\XX,\YY$). By \cref{lem:count of matching "blowups"} (applied with $\alpha=\eta$ and $W$ being the vertices appearing in $\XX$ or $\YY$), there are $\left(\Omega_{s}(1)\cdot  \frac{d^{s}}{(n/(\eta t))^{(s-4)}}\right)^{\eta ts/4}$ ways to do this. 
Putting together the above, if $n$ is large enough, for any fixed $\XX$, the number of tuples $(\XX,\YY,(M_{i,j})_{i,j\in[s/2]},\tilde\XX,\tilde \YY,(\tilde{M}_{i,j})_{i,j\in[s/2]})$ we get is at least

 $$\Omega_{s}(\eta)^{\eta ts^2/2}{\deg_\AA(\XX)}\cdot \left(\Omega_{s}(1)\cdot  \frac{d^{s}}{(n/(\eta t))^{(s-4)}}\right)^{\eta ts/4}.$$

Hence, the double counting of the tuples gives
$$\sum_{\XX\textnormal{ bad}}\deg_\AA(\XX)\cdot \Omega_{s}(\eta )^{\eta ts^2/2}\cdot \left(\frac{d^{s}}{(n/(\eta t))^{(s-4)}}\right)^{\eta ts/4}
\leq 2e(\AA)\cdot \left(O_s\left(\frac{n}{\eta^3 t}\right)\right)^{\eta ts/2},$$
and so, using $n/t=\Theta_s(d^{\frac{s}{s-2}}(\log n)^{-45})$,
\begin{align*}
    \sum_{\XX\textnormal{ bad}}\deg_\AA(\XX)&\leq  e(\AA)\cdot \left(O_s\left(\eta^{-(3s+2)}\frac{(n/t)^{s-2}}{d^{s} }\right)\right)^{\eta ts/4}\\
    &=e(\AA)\cdot \left(O_s\left(\mu^{-3(3s+2)}(\log n)^{-45(s-2)}\right)\right)^{\eta ts/4}.
\end{align*}
However, recalling $\mu=1/(\log n)^2$, we have 
$$\mu^{-3(3s+2)}(\log n)^{-45(s-2)}= (\log n)^{(18s+12)-(45s-90)}= (\log n)^{102-27s}=o(1).$$
Combining with $\eta t=\Omega_s((\log n)^{15})$, the claim follows.
    \renewcommand{\qedsymbol}{$\blacksquare$}
\end{proof}
    \renewcommand{\qedsymbol}{$\square$}

Recall that we took a random walk $\XX_{-N},\dots,\XX_N$ on $\AA$ starting from the stationary distribution. Using Claim~\ref{claim:few bad}, for each fixed $i\in[-N,N]$, $\XX_i$ is bad with probability at most $e^{-(\Omega_s(\log n))^{15}}$, which is at most $\frac{1}{6N}$ if $n$ is large enough, using $N=O_s(\mu^{-2})=O_s((\log n)^{4})$. Moreover, if $i\in [-(N-1),N-1]$ and we condition on $\XX_i$ being good, then (by the definition of $\XX_i$ being good and the observation that $\XX_{i-1}, \XX_{i+1}$ are independent $\AA$-neighbours of $\XX_i$ picked with probabilities proportional to $\AA$-weights), with probability at least $1-\frac{1}{6N}$, all but $s^3\eta t\leq \mu t/(6N)$ of the vertices in $\XX_{i}$ have degree at least $s$ in $\XX_{i-1}$ and $\XX_{i+1}$ combined. Hence, for each fixed $i\in [-(N-1),N-1]$, with probability at least $1-\frac{1}{3N}$, $\XX_i$ is good and has the property that all but $\mu t/(6N)$ of its vertices $\XX_{i}$ have degree at least $s$ in $\XX_{i-1}$ and $\XX_{i+1}$ combined. Therefore, with positive probability, for all $i\in[-(N-1),N-1]$ simultaneously, $\XX_i$ is good and has the property that all but $\mu t/(6N)$ of its vertices $\XX_{i}$ have degree at least $s$ in $\XX_{i-1}$ and $\XX_{i+1}$ combined; consider such an outcome.

Let $j\geq 1$ be minimal such that at most $(\mu/3) t$ of the vertices in $\XX_j$ do not appear in $\XX_0,\dots,\XX_{j-1}$, if there is such a $j$. Otherwise, let $j=N$. Similarly, let $\ell\leq -1$ be maximal such that at most $(\mu/3) t$ of the vertices in $\XX_\ell$ do not appear in $\XX_0,\XX_{-1},\dots,\XX_{\ell+1}$; if no such $\ell$ exists, let $\ell=-N$. We claim that the vertices spanned by $\XX_\ell,\dots,\XX_j$ satisfy the conditions. Let $U$ be the set of these vertices of $G$.

The number of vertices in $\XX_{\ell+1},\dots,\XX_{j-1}$ which have degree less than $s$ in $G[U]$ is at most $\mu t/(6N)\cdot 2N=(\mu/3) t$. We now count the number of vertices in $\XX_j,\XX_\ell$ with degree less than $s$ in $G[U]$. If $j<N$, then $\XX_j$ contains at most $(\mu/3)t$ additional vertices not appearing in $\XX_0,\dots,\XX_{j-1}$, thus, at most this many additional vertices of degree less than $s$. A similar statement holds if $\ell>-N$. Thus, if both $j<N$ and $\ell>-N$, then at most $(\mu/3)t+2\cdot (\mu /3)t= \mu t$ vertices have degree less than $s$, and we have at least $t$ vertices in total, giving the result. Therefore, we may assume that $j=N$ or $\ell=-N$; without loss of generality, $j=N$. But then the number of vertices appearing in $\XX_0,\dots,\XX_N$ is at least $N\cdot (\mu/3) t$. So $U$ contains at least $N(\mu/3)t$ vertices, out of which at most $(\mu/3)t+2(s/2)t=(s+\mu/3)t$ have degree less than $s$ in $G[U]$. Since we chose $N=\lceil (3s+1)/\mu ^2\rceil\geq 3(s+\mu/3)/\mu^2$, this gives that at most a proportion of $\frac{(s+\mu/3) t}{N(\mu /3) t}\leq \mu$ vertices have degree less than $s$, as required.
\end{proof}

\section{Odd $s$} \label{sec:odd}

In this section we prove Theorem \ref{thm:almost all degree s unpruned}. We first give an outline of the proof with more details than in the introduction. Recall that we want to prove that in every $n$-vertex graph $G$ with at least $n^{2-2/s+\eps}$ edges there exists a bounded-order subgraph in which almost all vertices have degree at least $s$. We will show that in fact each such graph $G$ contains one of three possible bounded-order subgraphs with few vertices of degree less than $s$. We will use three different embedding procedures, each of which will work under different assumptions about the counts of certain small subgraphs. Before we present these embedding procedures, let us introduce some notation. In what follows, let $s=2k+1$. Let $P_{\ell}$ denote the $\ell$-edge path. Moreover, let $H[t]$ denote the $t$-blowup of a graph $H$, that is, the graph obtained from $H$ by replacing each vertex with an independent set of size $t$ and each edge with a complete bipartite graph $K_{t,t}$. Finally, let $H\Box F$ denote the Cartesian product of graphs $H$ and $F$.

We now present our first embedding procedure, that will work if the subgraph counts are close to random. We will aim to find the graph $P_{\ell}[k]\Box K_2$ as a subgraph of $G$ (see Figure \ref{fig:embedding 1} for an illustration of this graph). Note that if $\ell$ is large, then almost all vertices of this graph have degree $2k+1$, as desired. The building blocks of this graph are copies of $K_{k,k}\Box K_2$, which is drawn on Figure \ref{fig:generalized cube}. We will attempt to build a copy of $P_{\ell}[k]\Box K_2$ by gluing together copies of $K_{k,k}\Box K_2$ in a chain along matchings of size $k$. Note that this is similar to how one could attempt to find a copy of $P_{\ell}[k]$ by gluing together copies of $K_{k,k}$ along independent sets of size $k$, but is considerably more difficult. In their proof of the average-degree analogue of our Theorem~\ref{thm:small odd subgraph}, Jiang and Newman \cite{jiang2017small} studied the graph $P_{\ell}[k]\Box K_2$ and showed that it can be found as a subgraph of $G$ if $G$ does not contain $K_{k+1,t}$ for some large constant $t$. Our starting point is a different argument that (under a technical condition about four-cycles) obtains a significantly stronger conclusion: informally, either $G$ contains a copy of $P_{\ell}[k]\Box K_2$, or it contains many copies of $K_{k,k+1}$ in which the part of size $k+1$ has large codegree -- we will call these \emph{fat} copies of $K_{k,k+1}$. (We refer the curious reader to Proposition \ref{prop:embedding 1} for the formal statement, which is more subtle than the simplified version we gave just now.)

\begin{figure}
	\centering
	\begin{tikzpicture}[scale=1.3]

	% Edges x_i y_i
	\draw[thick]
		(0,2)--(0,0)
		(1,2)--(1,0)
		(2,2)--(2,0)
		(3,2)--(3,0)
		(4,2)--(4,0)
		(5,2)--(5,0)
		(6,2)--(6,0);

	% Edges with |i-j|=1
	\draw[thick]
		(0,2)--(1,0) (1,2)--(0,0)
		(1,2)--(2,0) (2,2)--(1,0)
		(2,2)--(3,0) (3,2)--(2,0)
		(3,2)--(4,0) (4,2)--(3,0)
		(4,2)--(5,0) (5,2)--(4,0)
		(5,2)--(6,0) (6,2)--(5,0);

	% Edges with |i-j|=2
	\draw[thick]
		(0,2)--(2,0) (2,2)--(0,0)
		(1,2)--(3,0) (3,2)--(1,0)
		(2,2)--(4,0) (4,2)--(2,0)
		(3,2)--(5,0) (5,2)--(3,0)
		(4,2)--(6,0) (6,2)--(4,0);

	% Vertices
	\foreach \i in {0,...,6} {
		\draw[fill=black] (\i,2) circle (4pt);
		\draw[fill=black] (\i,0) circle (4pt);
	}

	% Labels
	\node at (0,2.4) {$x_1$};
	\node at (1,2.4) {$x_2$};
	\node at (2,2.4) {$x_3$};
	\node at (3,2.4) {$x_4$};
	\node at (4,2.4) {$x_5$};
	\node at (5,2.4) {$x_6$};
	\node at (6,2.4) {$x_7$};

	\node at (0,-0.4) {$y_1$};
	\node at (1,-0.4) {$y_2$};
	\node at (2,-0.4) {$y_3$};
	\node at (3,-0.4) {$y_4$};
	\node at (4,-0.4) {$y_5$};
	\node at (5,-0.4) {$y_6$};
	\node at (6,-0.4) {$y_7$};

	\end{tikzpicture}
	\caption{The graph $F_{2,7}$}
	\label{fig:F27}
\end{figure}
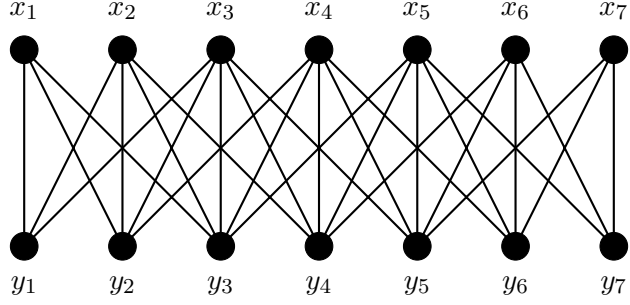

After this, we consider two cases. If $G$ contains few copies of $K_{k,k}$, then we can exploit the abundance of fat copies of $K_{k,k+1}$ to build a certain subgraph, depicted in Figure \ref{fig:F27}, that has few vertices of degree less than $s$. This is what we will call embedding procedure 3, and it is formalized by Proposition \ref{prop:embedding 3}. Finally, if $G$ contains many copies of $K_{k,k}$, then by a simple homomorphism inequality, it must contain many copies of $K_{k,k+1}$ as well. Let us address a very subtle point: we had already assumed that $G$ contains many fat copies of $K_{k,k+1}$, but here our conclusion is that $G$ has even more, but not necessarily fat copies of $K_{k,k+1}$. In particular, $G$ contains a lot more copies of $K_{k,k+1}$ than copies of $K_{1,k+1}$, which will imply that it has much more copies of $K_{k+1,k+1}$ than copies of $K_{1,k+1}$. Adapting an argument from~\cite{grzesik2022turan}, one can show that under this assumption $G$ contains a copy of $P_{\ell}[k+1]$ for some large constant~$\ell$. This is what we will call embedding procedure 2 (see Proposition \ref{prop:embedding 2}). Together, these propositions will cover all possible cases of the small subgraph counts and find a desired subgraph in each case.

\medskip

We now turn to the formal proof of Theorem \ref{thm:almost all degree s unpruned}.
The following lemma from \cite{chakraborti2024regular} will allow us to restrict our attention to regular host graphs.

\begin{lemma}[{\cite[Theorem 1.15]{chakraborti2024regular}}] \label{lemma:regular subgraph}
    For any positive real $\eps<1$, there exists some $\beta>0$ such that the following holds. Let $c>0$ and let $n$ be a positive integer that is sufficiently large as a function of $\eps$ and $c$. Let $G$ be a graph on $n$ vertices with $e(G)\geq cn^{1+\eps}$. Then $G$ contains a regular subgraph $H$ on $m\geq n^{\beta}$ vertices such that $e(H)\geq \beta c m^{1+\eps}$.
\end{lemma}

Using the simple fact that every graph contains a bipartite subgraph with at least half as many edges together with
Lemma \ref{lemma:regular subgraph}, Theorem \ref{thm:almost all degree s unpruned} reduces to the following statement.

\begin{proposition} \label{prop:almost all degree s}
    Let $s\geq 5$ be odd. Then for any $\eps,\mu>0$, there exist some $t=t(s,\mu)$ and $n_0=n_0(s,\eps,\mu)$ such that for each $n\geq n_0$, every $n$-vertex $d$-regular bipartite graph with $d\geq n^{1-2/s+\eps}$ contains a (non-empty) subgraph on at most $t$ vertices in which all but an at most $\mu$ proportion of vertices have degree at least $s$.
\end{proposition}

For the rest of this section, let $s=2k+1$, where $k\geq 2$ is an integer.

\begin{definition} \label{defn:fat}
    We call a copy of $K_{k,k+1}$ in an $n$-vertex graph $G$ \emph{fat} if the codegree of the part of size $k+1$ is at least $\log n$.
\end{definition}

We shall now state our three embedding results. 

\begin{proposition}[Embedding procedure 1] \label{prop:embedding 1}
    Let $\eps>0$, let $k\geq 2$ and $\ell$ be fixed positive integers and let $G$ be a $d$-regular graph on $n$ vertices, where $d\geq n^{1-2/(2k+1)+\eps}$. Let $\lambda=(nd^{k-3})^{1/(k-1)}\log n$. Let $h$ be the number of $4$-cycles in $G$ in which each diagonal pair of vertices has codegree at most $\lambda$. If $h=\Omega(d^4)$, then, for sufficiently large $n$, either $G$ contains the graph $P_{\ell}[k]\Box K_2$ or it contains at least $n^{-o(1)}h^{k^2}(nd)^{2k-2k^2}/\lambda^k$ fat copies of $K_{k,k+1}$.
\end{proposition}

Note that all but a $2/(\ell+1)$ proportion of the vertices of $P_{\ell}[k]\Box K_2$ have degree $2k+1=s$, so it is a suitable subgraph for Proposition \ref{prop:almost all degree s}, as long as $\ell$ is large enough.

\begin{proposition}[Embedding procedure 2] \label{prop:embedding 2}
    For any positive integers $k$ and $\ell$ and any integer $2\leq j\leq k+1$, there exists some $C$ such that the following holds. Let $G$ be a graph in which the number of copies of $K_{j,k+1}$ is greater than $C$ times the number of copies of $K_{1,k+1}$. Then $G$ contains $P_{\ell}[k+1]$.
\end{proposition}

Note that all but a $2/(\ell+1)$ proportion of the vertices of $P_{\ell}[k+1]$ have degree $2k+2>s$, so it is also a suitable subgraph for Proposition \ref{prop:almost all degree s}, as long as $\ell$ is large enough.

For our last embedding result, we need to make the following definition (see Figure \ref{fig:F27}).

\begin{definition}
    For positive integers $k$ and $\ell$, let $F_{k,\ell}$ be the graph with vertex set $\{x_1,\dots,x_{\ell},y_1,\dots,y_{\ell}\}$ and edge set $\{x_iy_j: |i-j|\leq k\}$.
\end{definition}

Observe that all but at most $2k/\ell$ proportion of the vertices in $F_{k,\ell}$ have degree $2k+1=s$.

\begin{proposition}[Embedding procedure 3] \label{prop:embedding 3}
    For any positive integers $k$ and $\ell$, there exists some $C$ such that the following holds for all sufficiently large $n$. Let $G$ be an $n$-vertex bipartite graph in which the number of fat copies of $K_{k,k+1}$ is greater than $C\log n$ times the number of copies of $K_{k,k}$. Then $G$ contains $F_{k,\ell}$.
\end{proposition}

We will now complete the proof of Proposition \ref{prop:almost all degree s}, assuming Propositions \ref{prop:embedding 1}, \ref{prop:embedding 2} and \ref{prop:embedding 3}.
We will first prove a useful corollary of Proposition \ref{prop:embedding 2}.

\begin{lemma} \label{lem:codegree bound}
    Let $k\geq 2$ and $\ell$ be fixed positive integers and let $G$ be a (non-empty) $d$-regular graph on $n$ vertices. Let $\lambda=(nd^{k-3})^{1/(k-1)}\log n$. If there exist $\Omega(d^4)$ many $4$-cycles in $G$ in which some diagonal pair of vertices has codegree greater than $\lambda$, then, for sufficiently large $n$, $G$ contains $P_{\ell}[k+1]$.
\end{lemma}

\begin{proof}
    Under the assumptions of the lemma, the number of copies of $K_{2,k+1}$ in $G$ is $\Omega(d^4\cdot \lambda^{k-1})=\omega(nd^{k+1})$ since any $4$-cycle in which some diagonal pair of vertices has codegree greater than $\lambda$ can be extended to $\Omega(\lambda^{k-1})$ copies of $K_{2,k+1}$. On the other hand, the number of copies of $K_{1,k+1}$ is at most $nd^{k+1}$. Hence, an application of Proposition \ref{prop:embedding 2} yields the desired conclusion.
\end{proof}

We recall the well-known result (see, e.g., \cite{erdHos1983supersaturated}) that $C_4$ supersaturates above the extremal number.

\begin{lemma} \label{lem:many 4cycles}
    Let $G$ be an $n$-vertex graph with average degree $d=\omega(n^{1/2})$. Then, for $n$ sufficiently large, $G$ contains $\Omega(d^4)$ copies of $C_4$.
\end{lemma}

We will now prove that either we can find a suitable subgraph or there are many copies of $K_{k,k}$ in our graph.

\begin{lemma} \label{lem:many k,k}
    Let $\eps>0$, let $k\geq 2$ and $\ell$ be fixed positive integers and let $G$ be a $d$-regular bipartite graph on $n$ vertices, where $d\geq n^{1-2/(2k+1)+\eps}$. Let $\lambda=(nd^{k-3})^{1/(k-1)}\log n$. Let $h$ be the number of $4$-cycles in $G$. Then, for sufficiently large $n$, either $G$ contains one of the graphs $P_{\ell}[k]\Box K_2$, $P_{\ell}[k+1]$ and $F_{k,\ell}$, or it contains at least $n^{-o(1)}h^{k^2}(nd)^{2k-2k^2}/\lambda^k$ copies of $K_{k,k}$.
\end{lemma}

\begin{proof}
    Let $h$ be the number of $4$-cycles in $G$. By Lemma \ref{lem:many 4cycles}, we have $h=\Omega(d^4)$. If at least $h/2$ four-cycles have a diagonal pair of vertices with codegree greater than $\lambda$, then Lemma \ref{lem:codegree bound} implies that $G$ contains $P_{\ell}[k+1]$. Else, it has at least $h/2$ four-cycles in which the diagonal pairs have codegree at most $\lambda$. Hence, by Proposition \ref{prop:embedding 1}, either $G$ contains $P_{\ell}[k]\Box K_2$ or it contains at least $n^{-o(1)}h^{k^2}(nd)^{2k-2k^2}/\lambda^k$ fat copies of $K_{k,k+1}$. In the latter case, Proposition \ref{prop:embedding 3} implies that either $G$ contains $F_{k,\ell}$, or it has at least $n^{-o(1)}h^{k^2}(nd)^{2k-2k^2}/\lambda^k$ copies of $K_{k,k}$.
\end{proof}

We are now nearly ready to prove Proposition \ref{prop:almost all degree s}. The proof is slightly different in the case $s=5$ (and is a bit simpler than the general case).

\begin{proposition} \label{prop:almost all degree 5}
    For any $\eps,\mu>0$, there exist some $t=t(\mu)$ and $n_0=n_0(\eps,\mu)$ such that for each $n\geq n_0$, every $n$-vertex $d$-regular bipartite graph with $d\geq n^{3/5+\eps}$ contains a subgraph on at most $t$ vertices in which all but an at most $\mu$ proportion of vertices have degree at least $5$.
\end{proposition}

\begin{proof}
    Given $\mu$, choose the smallest $\ell$ such that in each of the graphs $P_{\ell}[2]\Box K_2$, $P_{\ell}[3]$ and $F_{2,\ell}$, the proportion of vertices of degree less than $5$ is at most $\mu$. Let $h$ be the number of $4$-cycles in $G$ and let $\lambda=(n/d)\log n$, as in Lemma \ref{lem:many k,k} (for $k=2$). Then, by Lemma \ref{lem:many k,k}, assuming that $G$ does not contain a suitable subgraph, it must have at least $n^{-o(1)}h^4(nd)^{-4}/\lambda^2$ copies of $K_{2,2}$. Since $K_{2,2}$ is the same as $C_4$, it follows that $n^{-o(1)}h^4(nd)^{-4}/\lambda^2\leq h$, and so $h^3\leq n^{o(1)}(nd)^4\lambda^2\leq n^{o(1)}n^6 d^2$. On the other hand, by Lemma \ref{lem:many 4cycles}, we have $h=\Omega(d^4)$. Hence, we obtain $d^{12}\leq n^{o(1)}n^6 d^2$, which contradicts $d\geq n^{3/5+\eps}$.
\end{proof}

Turning to the case of general $s$, we now prove that if the number of copies of $K_{k,k}$ is large, then so is the number of copies of $K_{k,k+1}$. Here and below, $\hom(H,G)$ denotes the number of graph homomorphisms from $H$ to $G$.

\begin{lemma} \label{lem:simple hom inequality}
    We have $$\hom(K_{k,k+1},G)^{k-1}\cdot \hom(K_{k,1},G)\geq \hom(K_{k,k},G)^k.$$
\end{lemma}

\begin{proof}
    Observe that for any $1\leq i\leq k+1$, we have
    $$\hom(K_{k,i},G)=\sum_{f:[k]\rightarrow V(G)} |N(f([k]))|^i,$$
    where $N(S)$ denotes the common neighbourhood of the set $S$. By Hölder's inequality, we have 
    $$\sum_{f:[k]\rightarrow V(G)} |N(f([k]))|^{1/k}|N(f([k]))|^{(k^2-1)/k}\leq \left(\sum_{f:[k]\rightarrow V(G)} |N(f([k]))|\right)^{1/k}\left(\sum_{f:[k]\rightarrow V(G)} |N(f([k]))|^{k+1}\right)^{(k-1)/k},$$
    giving $\hom(K_{k,k},G)\leq \hom(K_{k,1},G)^{1/k}\cdot \hom(K_{k,k+1},G)^{(k-1)/k}$, as claimed.
    %Hence, the lemma follows from H\"older's inequality.
\end{proof}

We are now ready to prove Proposition \ref{prop:almost all degree s}.

\begin{proof}[Proof of Proposition \ref{prop:almost all degree s}]
    Since Proposition \ref{prop:almost all degree 5} deals with the case $s=5$, we may assume that $s=2k+1$ for some $k\geq 3$.

    Given $\mu$, choose the smallest $\ell$ such that in each of the graphs $P_{\ell}[k]\Box K_2$, $P_{\ell}[k+1]$ and $F_{k,\ell}$, the proportion of vertices of degree less than $s$ is at most $\mu$. Let $h$ be the number of $4$-cycles in $G$ and let $\lambda=(nd^{k-3})^{1/(k-1)}\log n$, as in Lemma \ref{lem:many k,k}. Then, by Lemma \ref{lem:many k,k}, assuming that $G$ does not contain a suitable subgraph, it must have at least $n^{-o(1)}h^{k^2}(nd)^{2k-2k^2}/\lambda^k$ copies of $K_{k,k}$.

    \medskip

    \noindent \textbf{Claim.} We have
    $$h^{k^2} (nd)^{2k-2k^2}/\lambda^k > n^{\Omega(1)}\cdot nd^{k+1-\frac{1}{k}}.$$

    \medskip

    \noindent \textbf{Proof of Claim.} Using $h=\Omega(d^4)$ and $\lambda=(nd^{k-3})^{1/(k-1)}\log n$, the desired inequality translates to
    $$n^{2k-2k^2-\frac{k}{k-1}}d^{2k+2k^2-\frac{k(k-3)}{k-1}}>n^{\Omega(1)}\cdot nd^{k+1-\frac{1}{k}},$$
    where the polylogarithmic term is absorbed by the $n^{\Omega(1)}$ term.
    Since $d\geq n^{1-\frac{2}{2k+1}}$, it suffices to prove that
    $$2k-2k^2-\frac{k}{k-1}+\left(1-\frac{2}{2k+1}\right)\left(2k+2k^2-\frac{k(k-3)}{k-1}\right)>1+\left(1-\frac{2}{2k+1}\right)\left(k+1-\frac{1}{k}\right).$$
    A simple computation simplifies this to $\frac{k^2-3k+1}{k(k-1)(2k+1)}>0$, which indeed holds for all $k\geq 3$. $\Box$

    \medskip

    Now by the claim, we have $\hom(K_{k,k},G)\geq n^{\Omega(1)}\cdot nd^{k+1-1/k}$ (absorbing the $n^{-o(1)}$ factor). Using Lemma \ref{lem:simple hom inequality} and the fact that $\hom(K_{k,1},G)=nd^k$, this implies that $\hom(K_{k,k+1},G)\geq n^{\Omega(1)}\cdot nd^{k+1}$, as well as that $\hom(K_{k,k+1},G)\geq n^{\Omega(1)}\cdot \hom(K_{k,k},G)$. The latter implies that at least half of the homomorphisms $K_{k,k+1}\rightarrow G$ are injective on the part of size $k+1$. (Indeed, any homomorphism $K_{k,k+1}\rightarrow G$ that is not injective on the part of size $k+1$ yields a homomorphism $K_{k,k}\rightarrow G$ with the same image, and each homomorphism $K_{k,k}\rightarrow G$ is obtained $O_k(1)$ times this way.) Hence, there are at least $n^{\Omega(1)}\cdot nd^{k+1}$ homomorphisms $K_{k,k+1}\rightarrow G$ that are injective on the part of size $k+1$. This implies that for some $1\leq j\leq k$, there are at least $n^{\Omega(1)}\cdot nd^{k+1}$ copies of $K_{j,k+1}$ in $G$. Since the number of copies of $K_{1,k+1}$ is at most $nd^{k+1}$, it follows from Proposition \ref{prop:embedding 2} that $G$ contains $P_{\ell}[k+1]$.
\end{proof}

\subsection{Embedding procedure 1}

\begin{figure}
	\centering
	\begin{tikzpicture}[scale=1.3]
	
	\draw[fill=black](0,0)circle(4pt);
    \draw[fill=black](0,1)circle(4pt);
    \draw[fill=black](0,2)circle(4pt);
    \draw[fill=black](0,3)circle(4pt);
    \draw[fill=black](0,4)circle(4pt);
    \draw[fill=black](0,5)circle(4pt);

	\draw[fill=black](2,0)circle(4pt);
    \draw[fill=black](2,1)circle(4pt);
    \draw[fill=black](2,2)circle(4pt);
    \draw[fill=black](2,3)circle(4pt);
    \draw[fill=black](2,4)circle(4pt);
    \draw[fill=black](2,5)circle(4pt);
    
    \draw[fill=black](4,0)circle(4pt);
    \draw[fill=black](4,1)circle(4pt);
    \draw[fill=black](4,2)circle(4pt);
    \draw[fill=black](4,3)circle(4pt);
    \draw[fill=black](4,4)circle(4pt);
    \draw[fill=black](4,5)circle(4pt);

    \draw[fill=black](6,0)circle(4pt);
    \draw[fill=black](6,1)circle(4pt);
    \draw[fill=black](6,2)circle(4pt);
    \draw[fill=black](6,3)circle(4pt);
    \draw[fill=black](6,4)circle(4pt);
    \draw[fill=black](6,5)circle(4pt);

    \draw[thick](0,0)--(0,1)(0,2)--(0,3)(0,4)--(0,5);
    \draw[thick](2,0)--(2,1)(2,2)--(2,3)(2,4)--(2,5);
    \draw[thick](4,0)--(4,1)(4,2)--(4,3)(4,4)--(4,5);
    \draw[thick](6,0)--(6,1)(6,2)--(6,3)(6,4)--(6,5);
    
    \draw[thick](0,0)--(2,0)(0,1)--(2,1)(0,2)--(2,2)(0,3)--(2,3)(0,4)--(2,4)(0,5)--(2,5);
    \draw[thick](0,0)--(2,2)(0,1)--(2,3)(0,2)--(2,4)(0,3)--(2,5)(0,4)--(2,0)(0,5)--(2,1);
    \draw[thick](0,0)--(2,4)(0,1)--(2,5)(0,2)--(2,0)(0,3)--(2,1)(0,4)--(2,2)(0,5)--(2,3);

    \draw[thick](4,0)--(2,0)(4,1)--(2,1)(4,2)--(2,2)(4,3)--(2,3)(4,4)--(2,4)(4,5)--(2,5);
    \draw[thick](4,0)--(2,2)(4,1)--(2,3)(4,2)--(2,4)(4,3)--(2,5)(4,4)--(2,0)(4,5)--(2,1);
    \draw[thick](4,0)--(2,4)(4,1)--(2,5)(4,2)--(2,0)(4,3)--(2,1)(4,4)--(2,2)(4,5)--(2,3);

    \draw[thick](4,0)--(6,0)(4,1)--(6,1)(4,2)--(6,2)(4,3)--(6,3)(4,4)--(6,4)(4,5)--(6,5);
    \draw[thick](4,0)--(6,2)(4,1)--(6,3)(4,2)--(6,4)(4,3)--(6,5)(4,4)--(6,0)(4,5)--(6,1);
    \draw[thick](4,0)--(6,4)(4,1)--(6,5)(4,2)--(6,0)(4,3)--(6,1)(4,4)--(6,2)(4,5)--(6,3);

    \node at (-2,2.5) {$\cdots$};
    \node at (8,2.5) {$\cdots$};

    \node at (-0.5,1) {$x^i_{3}$};
    \node at (-0.5,3) {$x^i_{2}$};
    \node at (-0.5,5) {$x^i_{1}$};
    \node at (-0.5,0) {$y^i_{3}$};
    \node at (-0.5,2) {$y^i_{2}$};
    \node at (-0.5,4) {$y^i_{1}$};

    \node at (2,5.4) {$x^{i+1}_{1}$};
    \node at (2,-0.4) {$y^{i+1}_{3}$};

    \node at (6.6,1) {$x^{i+3}_{3}$};
    \node at (6.6,3) {$x^{i+3}_{2}$};
    \node at (6.6,5) {$x^{i+3}_{1}$};
    \node at (6.6,0) {$y^{i+3}_{3}$};
    \node at (6.6,2) {$y^{i+3}_{2}$};
    \node at (6.6,4) {$y^{i+3}_{1}$};

    \node at (4,5.4) {$x^{i+2}_{1}$};
    \node at (4,-0.4) {$y^{i+2}_{3}$};

	\end{tikzpicture}		
	\caption{Embedding procedure 1 in the $s=7$ case}
	\label{fig:embedding 1}
\end{figure}
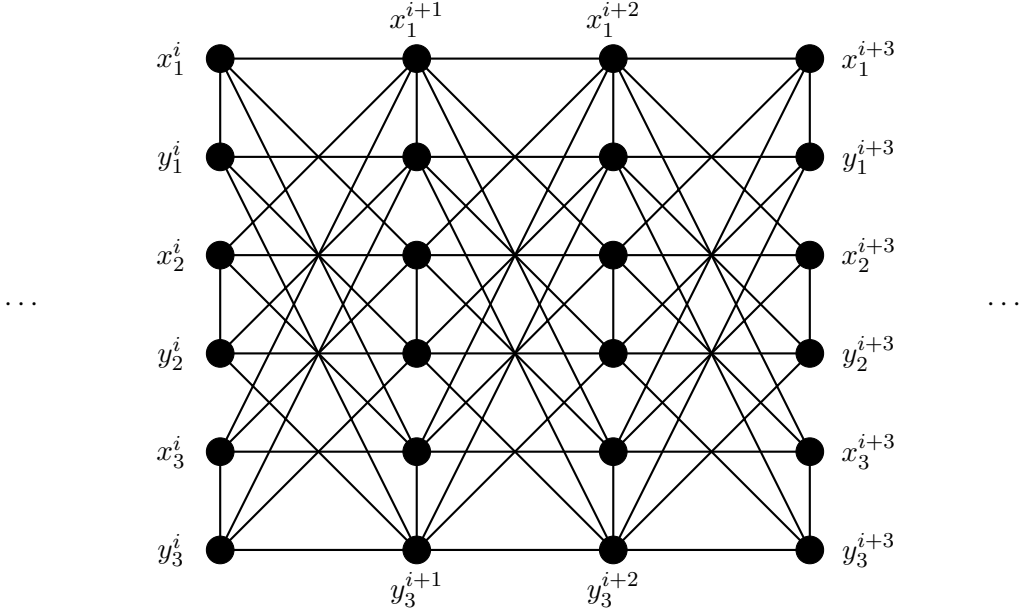

In this subsection, we define our first embedding procedure and prove Proposition \ref{prop:embedding 1}.
See Figure \ref{fig:embedding 1} for the graph we find in this case. The building blocks of this graph are copies of $K_{k,k}\Box K_2$, which is drawn on Figure \ref{fig:generalized cube}. Let us label the vertices of $K_{k,k}\Box K_2$ as in that figure. Moreover, let us label the vertices of $K_{k,1}\Box K_2$ as in Figure \ref{fig:star and edge}.

\begin{figure}
\centering
\begin{minipage}{0.48\textwidth}
\centering
\begin{tikzpicture}[scale=0.9]

\draw[fill=black](0,-2)circle(4pt);
    \draw[fill=black](0,-1)circle(4pt);
    \draw[fill=black](0,2)circle(4pt);
    \draw[fill=black](0,3)circle(4pt);
    \draw[fill=black](0,4)circle(4pt);
    \draw[fill=black](0,5)circle(4pt);

	\draw[fill=black](2,-2)circle(4pt);
    \draw[fill=black](2,-1)circle(4pt);
    \draw[fill=black](2,2)circle(4pt);
    \draw[fill=black](2,3)circle(4pt);
    \draw[fill=black](2,4)circle(4pt);
    \draw[fill=black](2,5)circle(4pt);
    
    \draw[fill=black](0,0.5)circle(1pt);
    \draw[fill=black](0,0.7)circle(1pt);
    \draw[fill=black](0,0.3)circle(1pt);

    \draw[fill=black](2,0.5)circle(1pt);
    \draw[fill=black](2,0.7)circle(1pt);
    \draw[fill=black](2,0.3)circle(1pt);

    \draw[thick](0,-2)--(0,-1)(0,2)--(0,3)(0,4)--(0,5);
    \draw[thick](2,-2)--(2,-1)(2,2)--(2,3)(2,4)--(2,5);
    
    \draw[thick](0,-2)--(2,-2)(0,-1)--(2,-1)(0,2)--(2,2)(0,3)--(2,3)(0,4)--(2,4)(0,5)--(2,5);
    \draw[thick](0,-2)--(2,2)(0,-1)--(2,3)(0,2)--(2,4)(0,3)--(2,5)(0,4)--(2,-2)(0,5)--(2,-1);
    \draw[thick](0,-2)--(2,4)(0,-1)--(2,5)(0,2)--(2,-2)(0,3)--(2,-1)(0,4)--(2,2)(0,5)--(2,3);

    \node at (-0.6,-1) {$x_k$};
    \node at (-0.6,3) {$x_2$};
    \node at (-0.6,5) {$x_1$};

    \node at (2.6,-2) {$z_k$};
    \node at (2.6,2) {$z_2$};
    \node at (2.6,4) {$z_1$};

    \node at (2.6,-1) {$w_k$};
    \node at (2.6,3) {$w_2$};
    \node at (2.6,5) {$w_1$};

    \node at (-0.6,-2) {$y_k$};
    \node at (-0.6,2) {$y_2$};
    \node at (-0.6,4) {$y_1$};

\end{tikzpicture}
\caption{A copy of $K_{k,k}\Box K_2$}
\label{fig:generalized cube}
\end{minipage}
\hfill
\begin{minipage}{0.48\textwidth}
\centering
\begin{tikzpicture}[scale=0.9]

\draw[fill=black](0,-2)circle(4pt);
    \draw[fill=black](0,-1)circle(4pt);
    \draw[fill=black](0,2)circle(4pt);
    \draw[fill=black](0,3)circle(4pt);
    \draw[fill=black](0,4)circle(4pt);
    \draw[fill=black](0,5)circle(4pt);

    \draw[fill=black](2,2)circle(4pt);
    \draw[fill=black](2,3)circle(4pt);
    
    \draw[fill=black](0,0.5)circle(1pt);
    \draw[fill=black](0,0.7)circle(1pt);
    \draw[fill=black](0,0.3)circle(1pt);

    \draw[thick](0,-2)--(0,-1)(0,2)--(0,3)(0,4)--(0,5);
    \draw[thick](2,2)--(2,3);
    
    \draw[thick](0,2)--(2,2)(0,3)--(2,3);
    \draw[thick](0,-2)--(2,2)(0,-1)--(2,3);
    \draw[thick](0,4)--(2,2)(0,5)--(2,3);

    \node at (-0.6,-1) {$x_k$};
    \node at (-0.6,3) {$x_2$};
    \node at (-0.6,5) {$x_1$};

    \node at (2.6,2) {$z$};

    \node at (2.6,3) {$w$};

    \node at (-0.6,-2) {$y_k$};
    \node at (-0.6,2) {$y_2$};
    \node at (-0.6,4) {$y_1$};

\end{tikzpicture}
\caption{A copy of $K_{k,1}\Box K_2$}
\label{fig:star and edge}
\end{minipage}
\end{figure}

The basic idea is to count certain homomorphic copies of $K_{k,k}\Box K_2$ and to glue them along matchings of size $k$.

\begin{definition} \label{def:light and weakly injective}
    Let us denote the vertices of $C_4$ as $x,y,z,w$ in the natural order. An injective homomorphism $f:V(C_4)\rightarrow V(G)$ is \emph{light} if the codegree of $f(x)$ and $f(z)$ is at most $\lambda$, and the codegree of $f(y)$ and $f(w)$ is at most $\lambda$, where $\lambda=(nd^{k-3})^{1/(k-1)}\log n$. 

    A homomorphism $f:V(K_{k,1}\Box K_2)\rightarrow V(G)$ is \emph{weakly injective} if $f|_{\{x_i,y_i,z,w\}}:V(C_4)\rightarrow V(G)$ is injective for all $i\in [k]$.

    A weakly injective homomorphism $f:V(K_{k,1}\Box K_2)\rightarrow V(G)$ is \emph{light} if the restricted (injective) homomorphism $f|_{\{x_i,y_i,z,w\}}:V(C_4)\rightarrow V(G)$ is light for all $i\in [k]$.

    A homomorphism $f:V(K_{k,k}\Box K_2)\rightarrow V(G)$ is \emph{weakly injective} if $f|_{\{x_i,y_i,z_j,w_j\}}:V(C_4)\rightarrow V(G)$ is injective for all $i,j\in [k]$.

    Finally, a weakly injective homomorphism $f:V(K_{k,k}\Box K_2)\rightarrow V(G)$ is \emph{light} if the restricted (injective) homomorphism $f|_{\{x_i,y_i,z_j,w_j\}}:V(C_4)\rightarrow V(G)$ is light for all $i,j\in [k]$.
\end{definition}

Given $H\in \{C_4,K_{k,1}\Box K_2,K_{k,k}\Box K_2\}$, let us write $\Hom^*(H,G)$ for the set of weakly injective and light graph homomorphisms $f:V(H)\rightarrow V(G)$, where we call a homomorphism $f:V(C_4)\rightarrow V(G)$ weakly injective if it is injective. Let us write $\hom^*(H,G)$ for the size of $\Hom^*(H,G)$.

We will use the following simple inequalities between these homomorphism counts.

\begin{lemma} \label{lem:hom inequalities}
    Let $G$ be a $d$-regular graph on $n$ vertices. Then
    $$\hom^*(K_{k,1}\Box K_2,G)\geq \frac{\hom^*(C_4,G)^k}{(nd)^{k-1}}$$ and
    $$\hom^*(K_{k,k}\Box K_2,G)\geq \frac{\hom^*(K_{k,1}\Box K_2,G)^k}{(nd)^{k^2-k}}.$$

    In particular, $$\hom^*(K_{k,k}\Box K_2,G)\geq \hom^*(C_4,G)^{k^2}(nd)^{2k-2k^2}.$$ 
\end{lemma}

\begin{proof}
    Given $u,v\in V(G)$ such that $uv\in E(G)$, let $\alpha_{u,v}$ be the number of weakly injective and light homomorphisms $f:V(C_4)\rightarrow V(G)$ such that $f(z)=u$ and $f(w)=v$. Observe that
    $$\hom^*(C_4,G)=\sum_{u,v\in V(G): uv\in E(G)} \alpha_{u,v}$$ and
    $$\hom^*(K_{k,1}\Box K_2,G)=\sum_{u,v\in V(G): uv\in E(G)} \alpha_{u,v}^k.$$
    Since the number of pairs $(u,v)\in V(G)^2$ with $uv\in E(G)$ is $nd$, we have, by Jensen's inequality, that
    $$\frac{1}{nd}\hom^*(K_{k,1}\Box K_2,G)\geq \left(\frac{1}{nd}\hom^*(C_4,G)\right)^k,$$
    which shows that the first inequality of the lemma holds.

    Now given $u_1,\dots,u_k,v_1,\dots,v_k\in V(G)$ such that $u_iv_i\in E(G)$ for each $i\in [k]$, let $\beta_{\textbf{u},\textbf{v}}$ be the number of weakly injective and light homomorphisms $f:V(K_{k,1}\Box K_2)\rightarrow V(G)$ such that $f(x_i)=u_i$ and $f(y_i)=v_i$ for all $i\in [k]$. Observe that
    $$\hom^*(K_{k,1}\Box K_2,G)=\sum_{\textbf{u}\in V(G)^k, \textbf{v}\in V(G)^k: u_iv_i\in E(G) \textrm { for all } i} \beta_{\textbf{u},\textbf{v}}$$ and
    $$\hom^*(K_{k,k}\Box K_2,G)=\sum_{\textbf{u}\in V(G)^k, \textbf{v}\in V(G)^k: u_iv_i\in E(G) \textrm { for all } i} \beta_{\textbf{u},\textbf{v}}^k.$$
    Since the number of pairs $(\textbf{u},\textbf{v})\in V(G)^k\times V(G)^k$ with $u_iv_i\in E(G)$ for all $i\in [k]$ is $(nd)^k$, we have, by Jensen's inequality, that
    $$\frac{1}{(nd)^k}\hom^*(K_{k,k}\Box K_2,G)\geq \left(\frac{1}{(nd)^k}\hom^*(K_{k,1}\Box K_2,G)\right)^k,$$
    which shows that the second inequality of the lemma holds.
    
    The third inequality follows by combining the first two inequalities.
\end{proof}

\begin{lemma} \label{lem:bounds for homs}
    Let $k\geq 2$ and let $G$ be a $d$-regular $n$-vertex graph with $d\geq n^{1-2/(2k+1)+\eps}$. If $\hom^*(C_4,G)=\Omega(d^4)$, then the following hold.
    \begin{enumerate}[(a)]
        \item \label{item:upper bound for k1} $\hom^*(K_{k,1}\Box K_2,G)=O(n^{-\eps}\hom^*(K_{k,k}\Box K_2,G))$. \label{item: (a)}
        \item $\hom^*(C_4,G)\lambda^{2k-2}=O(n^{-\eps}\hom^*(K_{k,k}\Box K_2,G))$. \label{item: (b)}
    \end{enumerate}
\end{lemma}

\begin{proof}

    Using Lemma \ref{lem:hom inequalities}, we have
    \begin{align*}
        \hom^*(K_{k,k}\Box K_2,G)
        &\geq \hom^*(K_{k,1}\Box K_2,G)\cdot \hom^*(K_{k,1}\Box K_2,G)^{k-1}\cdot (nd)^{k-k^2} \\
        &\geq \hom^*(K_{k,1}\Box K_2,G)\cdot \hom^*(C_4,G)^{k^2-k} (nd)^{-(k-1)^2}\cdot (nd)^{k-k^2} \\
        &\geq \hom^*(K_{k,1}\Box K_2,G)\cdot \Omega(d^{4(k^2-k)}) (nd)^{-(k-1)^2}\cdot (nd)^{k-k^2} \\
        &= \hom^*(K_{k,1}\Box K_2,G)\cdot \Omega(d^{2k^2-k-1}n^{-2k^2+3k-1}).
    \end{align*}
    Observe that $d^{2k^2-k-1}n^{-2k^2+3k-1}=(d^{2k+1}n^{-(2k-1)})^{k-1}=(dn^{-(1-2/(2k+1))})^{(2k+1)(k-1)}\geq n^{\eps}$, from which \ref{item: (a)} follows.

    Finally, $\hom^*(K_{k,k}\Box K_2,G)\geq \hom^*(C_4,G)^{k^2}(nd)^{2k-2k^2}$ and $\hom^*(C_4,G)=\Omega(d^4)$ imply that
    $\hom^*(K_{k,k}\Box K_2,G)\geq \hom^*(C_4,G)\cdot \Omega(d^{4k^2-4})\cdot (nd)^{2k-2k^2}=\hom^*(C_4,G)\cdot \Omega(d^{2k^2+2k-4} n^{2k-2k^2})$. On the other hand, $\lambda^{2k-2}=n^2 d^{2k-6}(\log n)^{2k-2}$. Hence, for \ref{item: (b)}, it is sufficient to prove that $d^{2k^2+2k-4} n^{2k-2k^2}\geq n^{\eps} n^2 d^{2k-6}(\log n)^{2k-2}$, for which it is enough to show that $d^{2k^2+2}\geq n^{2\eps} n^{2k^2-2k+2}$. This follows from $d\geq n^{1-2/(2k+1)+\eps}$.
\end{proof}

\begin{definition} \label{def:rich and poor}
    Let $G$ be an $n$-vertex graph and let $u_1,\dots,u_k,v_1,\dots,v_k\in V(G)$ such that there is at least one weakly injective and light homomorphism $f:V(K_{k,1}\Box K_2)\rightarrow V(G)$ with $f(x_i)=u_i$ and $f(y_i)=v_i$ for all $i\in [k]$. We call $(\textbf{u},\textbf{v})$ \emph{poor} if there exists some $q\in V(G)$ such that at least $1/\log n$ proportion of all such homomorphisms have $f(z)=q$ or $f(w)=q$, and we call $(\textbf{u},\textbf{v})$ \emph{rich} otherwise.

    \sloppy We call a weakly injective and light homomorphism $f:V(K_{k,k}\Box K_2)\rightarrow V(G)$ \emph{poor} if the vector $(f(x_1),\dots,f(x_k),f(y_1),\dots,f(y_k))$ is poor, and we call $f$ \emph{rich} otherwise.
\end{definition}

The relevance of this definition lies in the following result, showing that if the proportion of poor homomorphisms is small, then we can find the graph depicted in Figure \ref{fig:embedding 1}, which satisfies the conclusion of Proposition \ref{prop:almost all degree s}.

\begin{lemma} \label{lem:if few poor}
    For any fixed positive integers $k\geq 2$ and $\ell$, the following holds for all sufficiently large $n$. Suppose that the proportion of poor homomorphisms amongst all weakly injective and light homomorphisms $f:V(K_{k,k}\Box K_2)\rightarrow V(G)$ is at most $1/\log n$. Then $G$ contains a (genuine) copy of $P_\ell[k]\Box K_2$.
\end{lemma}

Let us label the vertices of $P_{\ell}[k]\Box K_2$ as in Figure \ref{fig:embedding 1}.

\begin{proof}
    Define an auxiliary graph $\mathcal{G}$ as follows. Let $V(\mathcal{G})$ consist of $(2k)$-tuples $(u_1,\dots,u_k,v_1,\dots,v_k)\in V(G)^{2k}$ such that $u_iv_i\in E(G)$ holds for all $i\in [k]$. Moreover, take an edge in $\GG$ between $(u_1,\dots,u_k,v_1,\dots,v_k)\in V(G)^{2k}$ and $(u'_1,\dots,u'_k,v'_1,\dots,v'_k)\in V(G)^{2k}$ if and only if the map $f:V(K_{k,k}\Box K_2)\rightarrow V(G)$ defined by $f(x_i)=u_i$, $f(y_i)=v_i$, $f(z_i)=v'_i$ and $f(w_i)=u'_i$ is a weakly injective and light graph homomorphism. Note that this definition is invariant under swapping $(u_1,\dots,u_k,v_1,\dots,v_k)$ and $(u'_1,\dots,u'_k,v'_1,\dots,v'_k)$, and moreover no vertex of $\GG$ is adjacent to itself, so $\GG$ is an undirected simple graph.

    Consider a random walk $(\vu^i,\vv^i)_{0\leq i\leq \ell}\in V(\GG)$, where $(\vu^0,\vv^0)$ is chosen randomly according to the stationary distribution on $\GG$ and $(\vu^{i+1},\vv^{i+1})$ is a uniformly random neighbour of $(\vu^i,\vv^i)$. Note that for each $0\leq i\leq \ell-1$ the edge between $(\vu^i,\vv^i)\in V(\GG)$ and $(\vu^{i+1},\vv^{i+1})\in V(\GG)$ is a uniformly random edge of $\GG$, with each endpoint of the edge being equally likely to be $(\vu^i,\vv^i)$ or $(\vu^{i+1},\vv^{i+1})$. Moreover, observe that the edges of $\GG$ are in one-to-two correspondence with the set of weakly injective and light homomorphisms $f:V(K_{k,k}\Box K_2)\rightarrow V(G)$ (swapping $f(x_i)$ with $f(w_i)$ and $f(y_i)$ with $f(z_i)$ leads to the same edge). Hence, for each $0\leq i\leq \ell-1$, the distribution of $(\vu^i,\vv^i)$ is the same as the distribution of $(f(x_1),\dots,f(x_k),f(y_1),\dots,f(y_k))$ for a uniformly random weakly injective and light homomorphism $f:V(K_{k,k}\Box K_2)\rightarrow V(G)$. Hence, by the assumption of the lemma and Definition \ref{def:rich and poor}, the probability that $(\vu^i,\vv^i)$ is poor is at most $1/\log n$.

    Notice that any outcome of $(\vu^0,\vv^0),\dots,(\vu^{\ell},\vv^{\ell})$ defines a random homomorphism $g:V(P_{\ell}[k]\Box K_2)\rightarrow V(G)$ by letting $g(\mathbf{x}^i,\mathbf{y}^i)=(\vu^i,\vv^i)$ for each $i$. However, $g$ may not be an injective homomorphism. We will now upper bound the probability that two vertices of $P_{\ell}[k]\Box K_2$ are mapped to the same vertex in~$G$.

    \medskip

    \noindent \textbf{Claim.} For any distinct vertices $a,b\in V(P_{\ell}[k]\Box K_2)$, we have $\mathbb{P}(g(a)=g(b))=o(1)$.

    \medskip

    \noindent \textbf{Proof of Claim.} Let $0\leq i,j\leq \ell$ satisfy that $a\in \{x^i_p: p\in [k]\}\cup \{y^i_p: p\in [k]\}$ and $b\in \{x^j_p: p\in [k]\}\cup \{y^j_p: p\in [k]\}$.
    
    Suppose first that $i<j$. Recall that the probability that $(\vu^{j-1},\vv^{j-1})$ is poor is at most $1/\log n$. Now condition on an outcome of $(\vu^0,\vv^0)$, $(\vu^1,\vv^1)$, \dots, $(\vu^{j-1},\vv^{j-1})$. Note that $(\vu^j,\vv^j)$ is a uniformly random neighbour of $(\vu^{j-1},\vv^{j-1})$.
    In fact, we can sample $(\vu^j,\vv^j)$ as follows. Let $S$ be the set of pairs $(w,z)\in V(G)^2$ such that $wz\in E(G)$ and $(u^{j-1}_p,v^{j-1}_p,z,w)$ is a light (genuine) $4$-cycle for all $p\in [k]$. Now for each $p\in [k]$, sample $(u^j_p,v^j_p)$ from $S$ independently and uniformly at random. It is easy to see that this way we obtain a uniformly random neighbour $(\vu^j,\vv^j)$ of $(\vu^{j-1},\vv^{j-1})$. Now, if $(\vu^{j-1},\vv^{j-1})$ is rich, by definition, the proportion of $(w,z)\in S$ with $w=g(a)$ or $z=g(a)$ is less than $1/\log n$. In particular, the probability (conditioned on $(\vu^0,\vv^0)$, $(\vu^1,\vv^1)$, \dots, $(\vu^{j-1},\vv^{j-1})$ such that $(\vu^{j-1},\vv^{j-1})$ is rich) that $g(b)=g(a)$ is less than $1/\log n$. Combining this with the probability that $(\vu^{j-1},\vv^{j-1})$ is poor, we obtain that
    $\mathbb{P}(g(a)=g(b))\leq 2/\log n$.

    The case where $i=j$ is very similar. By the reversibility of the random walk, we may assume that $j\geq \ell/2$, so in particular $j\neq 0$. Let $a\in \{x^j_p,y^j_p\}$ and $b\in \{x^j_q,y^j_q\}$. If $p=q$, then $a$ and $b$ are neighbours, so, as $g$ is a graph homomorphism, $\mathbb{P}(g(a)=g(b))=0$. Else, conditional on $(\vu^0,\vv^0)$, $(\vu^1,\vv^1)$, \dots, $(\vu^{j-1},\vv^{j-1})$, the pairs $(u^j_p,v^j_p)$ and $(u^j_q,v^j_q)$ are independent samples from $S$. Hence, if $(\vu^{j-1},\vv^{j-1})$ is rich, then once again the conditional probability that $g(a)=g(b)$ is at most $1/\log n$. Combining this with the probability that $(\vu^{j-1},\vv^{j-1})$ is poor, we once again obtain $\mathbb{P}(g(a)=g(b))\leq 2/\log n$. $\Box$

    \medskip

    Now by the claim and the union bound, the probability that $g$ is injective is positive (for large $n$), and hence $G$ contains a genuine copy of $P_{\ell}[k]\Box K_2$.
\end{proof}

Lemma \ref{lem:if few poor} shows that we can find a suitable subgraph if the proportion of poor homomorphisms amongst all weakly injective and light homomorphisms $f:V(K_{k,k}\Box K_2)\rightarrow V(G)$ is small. On the other hand, if the proportion of such homomorphisms is large, we can exploit it using the following lemma.

\begin{lemma} \label{lem:only w collapses}
    Let $k\geq 2$ and let $G$ be a $d$-regular $n$-vertex graph with $d\geq n^{1-2/(2k+1)+\eps}$. Assume that $\hom^*(C_4,G)=\Omega(d^4)$ and that the proportion of poor homomorphisms amongst all weakly injective and light homomorphisms $f:V(K_{k,k}\Box K_2)\rightarrow V(G)$ is greater than $1/\log n$. Then an $n^{-o(1)}$ proportion of all $f\in \Hom^*(K_{k,k}\Box K_2,G)$ satisfies that
    \begin{itemize}
        \item $f(w_i)=f(w_j)$ for all $i,j\in [k]$,
        \item $f(y_i)\neq f(y_j)$ for all $i\neq j\in [k]$ and
        \item $f(z_i)\neq f(z_j)$ for all $i\neq j\in [k]$.
    \end{itemize}
\end{lemma}

We remark that we could in addition require $f(x_i)\neq f(x_j)$ for all $i\neq j\in [k]$, but this will not be needed.

We refer the reader to Figure \ref{fig:only w collapses} for an illustration of the homomorphisms counted by this lemma. The important thing here is that with the vertices $x_i$ removed, we obtain a copy of $K_{k,k+1}$ (see Figure~\ref{fig:only w collapses and x removed}).

\begin{figure}[h]
\centering
\begin{minipage}{0.48\textwidth}
\centering
\begin{tikzpicture}[scale=0.9]

% Left column
\draw[fill=black](0,-2)circle(4pt);
\draw[fill=black](0,-1)circle(4pt);
\draw[fill=black](0,2)circle(4pt);
\draw[fill=black](0,3)circle(4pt);
\draw[fill=black](0,4)circle(4pt);
\draw[fill=black](0,5)circle(4pt);

% Right column
\draw[fill=black](3,-2)circle(4pt);
\draw[fill=black](3,2)circle(4pt);
\draw[fill=black](3,4)circle(4pt);

% w vertex
\draw[fill=black](3,5)circle(4pt);

% dots
\draw[fill=black](0,0.5)circle(1pt);
\draw[fill=black](0,0.7)circle(1pt);
\draw[fill=black](0,0.3)circle(1pt);
\draw[fill=black](3,0.5)circle(1pt);
\draw[fill=black](3,0.7)circle(1pt);
\draw[fill=black](3,0.3)circle(1pt);

% vertical edges
\draw[thick](0,-2)--(0,-1);
\draw[thick](0,2)--(0,3);
\draw[thick](0,4)--(0,5);

% y-z edges
\draw[thick](0,-2)--(3,-2);
\draw[thick](0,2)--(3,2);
\draw[thick](0,4)--(3,4);
\draw[thick](0,-2)--(3,2);
\draw[thick](0,2)--(3,4);
\draw[thick](0,4)--(3,-2);
\draw[thick](0,-2)--(3,4);
\draw[thick](0,2)--(3,-2);
\draw[thick](0,4)--(3,2);

% x-w edges
\draw[thick](0,5)--(3,5);
\draw[thick](0,3)--(3,5);
\draw[thick](0,-1)--(3,5);

% curved w-z edges
\draw[thick](3,5) to[bend left=20] (3,4);
\draw[thick](3,5) to[bend left=40] (3,2);
\draw[thick](3,5) to[bend left=60] (3,-2);

% labels
\node at (-0.6,-1) {$x_k$};
\node at (-0.6,3) {$x_2$};
\node at (-0.6,5) {$x_1$};

\node at (-0.6,-2) {$y_k$};
\node at (-0.6,2) {$y_2$};
\node at (-0.6,4) {$y_1$};

\node at (3,-2.5) {$z_k$};
\node at (3,1.5) {$z_2$};
\node at (3,3.5) {$z_1$};

\node at (4.6,5) {$w_1=\dots=w_k$};

\end{tikzpicture}
\caption{The homomorphisms counted by Lemma \ref{lem:only w collapses}}

\label{fig:only w collapses}
\end{minipage}
\hfill
\begin{minipage}{0.48\textwidth}
\centering
\begin{tikzpicture}[scale=0.9]

% y vertices
\draw[fill=black](0,-2)circle(4pt);
\draw[fill=black](0,2)circle(4pt);
\draw[fill=black](0,4)circle(4pt);

% z vertices
\draw[fill=black](3,-2)circle(4pt);
\draw[fill=black](3,2)circle(4pt);
\draw[fill=black](3,4)circle(4pt);

% w vertex
\draw[fill=black](3,5)circle(4pt);

% dots
\draw[fill=black](0,0.5)circle(1pt);
\draw[fill=black](0,0.7)circle(1pt);
\draw[fill=black](0,0.3)circle(1pt);
\draw[fill=black](3,0.5)circle(1pt);
\draw[fill=black](3,0.7)circle(1pt);
\draw[fill=black](3,0.3)circle(1pt);

% y-z edges
\draw[thick](0,-2)--(3,-2);
\draw[thick](0,2)--(3,2);
\draw[thick](0,4)--(3,4);
\draw[thick](0,-2)--(3,2);
\draw[thick](0,2)--(3,4);
\draw[thick](0,4)--(3,-2);
\draw[thick](0,-2)--(3,4);
\draw[thick](0,2)--(3,-2);
\draw[thick](0,4)--(3,2);

% curved w-z edges
\draw[thick](3,5) to[bend left=20] (3,4);
\draw[thick](3,5) to[bend left=40] (3,2);
\draw[thick](3,5) to[bend left=60] (3,-2);

% labels
\node at (-0.6,-2) {$y_k$};
\node at (-0.6,2) {$y_2$};
\node at (-0.6,4) {$y_1$};

\node at (3,-2.5) {$z_k$};
\node at (3,1.5) {$z_2$};
\node at (3,3.5) {$z_1$};

\node at (4.6,5) {$w_1=\dots=w_k$};

\end{tikzpicture}
\caption{Removing the vertices $x_i$ gives a copy of $K_{k,k+1}$}

\label{fig:only w collapses and x removed}
\end{minipage}

\end{figure}

We shall prove Lemma \ref{lem:only w collapses} in two steps. First, we prove that an $n^{-o(1)}$ proportion of all homomorphisms $f$ have $f(w_i)=f(w_j)$ for all $i,j\in [k]$.

\begin{lemma} \label{lem:w does collapse}
    Let $k\geq 2$ and let $G$ be an $n$-vertex graph. Assume that the proportion of poor homomorphisms amongst all weakly injective and light homomorphisms $f:V(K_{k,k}\Box K_2)\rightarrow V(G)$ is greater than $1/\log n$. Then an $n^{-o(1)}$ proportion of all $f\in \Hom^*(K_{k,k}\Box K_2,G)$ satisfies that $f(w_i)=f(w_j)$ for all $i,j\in [k]$.
\end{lemma}

\begin{proof}
    Let $u_1,\dots,u_k,v_1,\dots,v_k\in V(G)$ such that there is at least one weakly injective and light homomorphism $f:V(K_{k,1}\Box K_2)\rightarrow V(G)$ with $f(x_i)=u_i$ and $f(y_i)=v_i$ for all $i\in [k]$. Let us call $(\textbf{u},\textbf{v})$ \emph{$w$-poor} if there exists some $q\in V(G)$ such that at least $1/(2\log n)$ proportion of all such homomorphisms have $f(w)=q$, and let us call $(\textbf{u},\textbf{v})$ \emph{$z$-poor} if there exists some $q\in V(G)$ such that at least $1/(2\log n)$ proportion of all such homomorphisms have $f(z)=q$. Let us call a weakly injective and light homomorphism $f:V(K_{k,k}\Box K_2)\rightarrow V(G)$ \emph{$w$-poor} (respectively, $z$-poor) if the vector $(f(x_1),\dots,f(x_k),f(y_1),\dots,f(y_k))$ is $w$-poor (respectively, $z$-poor). Note that every poor pair $(\textbf{u},\textbf{v})$ is $w$-poor or $z$-poor, and hence every poor weakly injective and light homomorphism $f:V(K_{k,k}\Box K_2)\rightarrow V(G)$ is $w$-poor or $z$-poor. By symmetry, the number of $w$-poor homomorphisms is equal to the number of $z$-poor homomorphisms, so by the assumption of the lemma, the proportion of $w$-poor homomorphisms among all weakly injective and light homomorphisms $f:V(K_{k,k}\Box K_2)\rightarrow V(G)$ is greater than $1/(2\log n)$.

    Given $u_1,\dots,u_k,v_1,\dots,v_k\in V(G)$ such that $u_iv_i\in E(G)$ for each $i\in [k]$, let $\beta_{\textbf{u},\textbf{v}}$ be the number of weakly injective and light homomorphisms $g:V(K_{k,1}\Box K_2)\rightarrow V(G)$ such that $g(x_i)=u_i$ and $g(y_i)=v_i$ for all $i\in [k]$. Observe that the number of $w$-poor $f\in \Hom^*(K_{k,k}\Box K_2,G)$ is
    $$\sum_{\textbf{u}\in V(G)^k, \textbf{v}\in V(G)^k: (\mathbf{u},\mathbf{v}) \text{ is } w\text{-poor}} \beta_{\textbf{u},\textbf{v}}^k.$$
    On the other hand, the number of $f\in \Hom^*(K_{k,k}\Box K_2,G)$ with $f(w_i)=f(w_j)$ for all $i,j\in [k]$ is at least $$\sum_{\textbf{u}\in V(G)^k, \textbf{v}\in V(G)^k: (\mathbf{u},\mathbf{v}) \text{ is } w\text{-poor}} (2\log n)^{-k}\beta_{\textbf{u},\textbf{v}}^k.$$
    (Indeed, if $(\textbf{u},\textbf{v})$ is $w$-poor, then there are at least $(\beta_{\textbf{u},\textbf{v}}/(2\log n))^k$ such homomorphisms $f$ with $f((\textbf{x},\textbf{y}))=(\textbf{u},\textbf{v})$ and $f(w_i)=q$ for all $i$, where $q$ is the vertex witnessing that $(\textbf{u},\textbf{v})$ is $w$-poor.)

    It follows that the number of $f\in \Hom^*(K_{k,k}\Box K_2,G)$ with $f(w_i)=f(w_j)$ for all $i,j\in [k]$ is at least $n^{-o(1)}$ times the number of $w$-poor $f\in \Hom^*(K_{k,k}\Box K_2,G)$. Since the latter is at least $\hom^*(K_{k,k}\Box K_2,G)/(2\log n)$, the lemma follows.
\end{proof}

Given Lemma \ref{lem:w does collapse}, Lemma \ref{lem:only w collapses} follows from the following result (using symmetry).

\begin{lemma} \label{lem:few extra degeneracies}
    Let $k\geq 2$ and let $G$ be a $d$-regular $n$-vertex graph with $d\geq n^{1-2/(2k+1)+\eps}$. Assume that $\hom^*(C_4,G)=\Omega(d^4)$. Let $\mathcal{F}$ be the set of all $f\in \Hom^*(K_{k,k}\Box K_2,G)$ satisfying that $f(w_i)=f(w_j)$ for all $i,j\in [k]$.
    \begin{enumerate}[(i)]
        \item The number of $f\in \mathcal{F}$ with $f(y_1)=f(y_2)$ is at most $O(n^{-\eps/(k-1)} \hom^*(K_{k,k}\Box K_2,G))$. \label{item: (i)}
        \item The number of $f\in \mathcal{F}$ with $f(z_1)=f(z_2)$ is at most $O(n^{-\eps/(k-1)} \hom^*(K_{k,k}\Box K_2,G))$. \label{item: (ii)}
    \end{enumerate}
\end{lemma}

We prepare the proof of Lemma \ref{lem:few extra degeneracies} with the following two simple lemmas.

\begin{lemma} \label{lem:all equal and two equal}
    Let $X_1,\dots,X_k$ be independent and identically distributed random variables, taking values in a finite set. Then
    $\mathbb{P}(X_1=\dots=X_k)\geq \mathbb{P}(X_1=X_2)^{k-1}$.
\end{lemma}

\begin{proof}
    Let $S$ be the set of values taken by the random variables $X_i$, and let $p_s=\mathbb{P}(X_i=s)$. Then
    $$\mathbb{P}(X_1=\dots=X_k)=\sum_{s\in S} p_s^{k}$$
    and
    $$\mathbb{P}(X_1=X_2)=\sum_{s\in S} p_s^2.$$
    By H\"older's inequality, we have
    $$(\sum_{s\in S} p_s^k)^{\frac{1}{k-1}}\cdot (\sum_{s\in S} p_s)^{\frac{k-2}{k-1}}\geq \sum_{s\in S} p_s^2,$$
    so, noting that $\sum_{s\in S} p_s=1$, the lemma follows.
\end{proof}

\begin{lemma} \label{lem:all equal and two equal general}
    Let $X_1,\dots,X_k$ and $Q$ be random variables, taking values in a finite set. Assume that for any $q$ in the support of $Q$, we have that the conditional random variables $(X_i \mid Q=q)_{i\in [k]}$ are independent and identically distributed. Then $\mathbb{P}(X_1=\dots=X_k)\geq \mathbb{P}(X_1=X_2)^{k-1}$.
\end{lemma}

\begin{proof}
    We have
    \begin{align*}
        \mathbb{P}(X_1=\dots=X_k)
        &=\sum_q \mathbb{P}(Q=q)\mathbb{P}(X_1=\dots=X_k \mid Q=q) \\
        &\geq \sum_q \mathbb{P}(Q=q)\mathbb{P}(X_1=X_2 \mid Q=q)^{k-1} \\
        &\geq \left( \sum_q \mathbb{P}(Q=q)\mathbb{P}(X_1=X_2 \mid Q=q)\right)^{k-1} \\
        &=\mathbb{P}(X_1=X_2)^{k-1},
    \end{align*}
    where the first inequality follows from Lemma \ref{lem:all equal and two equal}, and the second inequality follows from Jensen's inequality.
\end{proof}

We are now ready to prove Lemma \ref{lem:few extra degeneracies}.

\begin{proof}[Proof of Lemma \ref{lem:few extra degeneracies}]
    We start with \ref{item: (i)}. Let $F\in \mathcal{F}$ be chosen uniformly at random. Note that the random variables $X_i=F(y_i)$ satisfy the conditions of Lemma \ref{lem:all equal and two equal general} with $Q=(F(z_1),\dots,F(z_k),F(w_1),\dots,F(w_k))$. Hence, by that lemma, we have $$\mathbb{P}(F(y_1)=F(y_2))^{k-1}\leq \mathbb{P}(F(y_1)=\dots=F(y_k)).$$
    Let $N$ be the number of homomorphisms $f\in \mathcal{F}$ with $f(y_1)=f(y_2)$, and let $N'$ be the number of homomorphisms $f\in \mathcal{F}$ with $f(y_1)=\dots=f(y_k)$. The above inequality gives $N^{k-1}\leq |\mathcal{F}|^{k-2} N'$.

    Now observe that $N'\leq \hom^*(C_4,G)\lambda^{2k-2}$, or in other words, the number of ways to choose $f\in \Hom^*(K_{k,k}\Box K_2,G)$ with $f(w_1)=\dots=f(w_k)$ and $f(y_1)=\dots=f(y_k)$ is at most $\hom^*(C_4,G)\lambda^{2k-2}$. Indeed, there are at most $\hom^*(C_4,G)$ ways to choose $(f(x_1),f(y_1),f(z_1),f(w_1))$, and any such choice determines $f(y_i)$ and $f(w_i)$ for all $i$. Now $f(x_i)$ and $f(z_i)$ must be common neighbours of $f(y_i)$ and $f(w_i)$, and the codegree of $f(y_i)$ and $f(w_i)$ is at most $\lambda$, since $f|_{\{x_i,y_i,z_i,w_i\}}:V(C_4)\rightarrow V(G)$ is light. Therefore, there are at most $\lambda$ choices for each $f(x_i)$ and $f(z_i)$ ($2\leq i\leq k$), so we indeed have $N'\leq \hom^*(C_4,G)\lambda^{2k-2}$.
    By Lemma \ref{lem:bounds for homs} \ref{item: (b)}, it follows that $N'=O(n^{-\eps}\hom^*(K_{k,k}\Box K_2,G))$, so $N^{k-1}\leq |\mathcal{F}|^{k-2}N'\leq O(|\mathcal{F}|^{k-2}n^{-\eps}\hom^*(K_{k,k}\Box K_2,G))\leq O(n^{-\eps} \hom^*(K_{k,k}\Box K_2,G)^{k-1})$. This proves~\ref{item: (i)}.

    \sloppy We now turn to \ref{item: (ii)}. Once again, let $F\in \mathcal{F}$ be chosen uniformly at random. Note that the random variables $X_i=F(z_i)$ satisfy the conditions of Lemma \ref{lem:all equal and two equal general} with $Q=(F(x_1),\dots,F(x_k),F(y_1),\dots,F(y_k),F(w_1),\dots,F(w_k))$.
    Hence, by that lemma, we have $$\mathbb{P}(F(z_1)=F(z_2))^{k-1}\leq \mathbb{P}(F(z_1)=\dots=F(z_k)).$$
    Let $M$ be the number of homomorphisms $f\in \mathcal{F}$ with $f(z_1)=f(z_2)$, and let $M'$ be the number of homomorphisms $f\in \mathcal{F}$ with $f(z_1)=\dots=f(z_k)$. The above inequality gives $M^{k-1}\leq |\mathcal{F}|^{k-2} M'$.

    Now observe that $M'\leq \hom^*(K_{k,1}\Box K_2,G)$. Indeed, any homomorphism $f\in \Hom^*(K_{k,k}\Box K_2,G)$ with $f(w_1)=\dots=f(w_k)$ and $f(z_1)=\dots=f(z_k)$ can be viewed as a homomorphism in $\Hom^*(K_{k,1}\Box K_2,G)$. Hence, by Lemma \ref{lem:bounds for homs} \ref{item: (a)}, we have $M'=O(n^{-\eps} \hom^*(K_{k,k}\Box K_2,G))$, so
    $M^{k-1}\leq |\mathcal{F}|^{k-2} M'\leq O(|\mathcal{F}|^{k-2} n^{-\eps}\hom^*(K_{k,k}\Box K_2,G))\leq O(n^{-\eps} \hom^*(K_{k,k}\Box K_2,G)^{k-1})$. This proves~\ref{item: (ii)}.
\end{proof}

We are now ready to prove Proposition \ref{prop:embedding 1}.

\begin{proof}[Proof of Proposition \ref{prop:embedding 1}]
    Note that $\hom^*(C_4,G)=\Theta(h)$. Hence, by the assumption of the proposition, we have $\hom^*(C_4,G)=\Omega(d^4)$. If the proportion of poor homomorphisms amongst all weakly injective and light homomorphisms $f:V(K_{k,k}\Box K_2)\rightarrow V(G)$ is at most $1/\log n$, then Lemma~\ref{lem:if few poor} implies that $G$ contains a copy of $P_\ell[k]\Box K_2$ and we are done. Else, Lemma \ref{lem:only w collapses} implies that an $n^{-o(1)}$ proportion of all $f\in \Hom^*(K_{k,k}\Box K_2,G)$ satisfies that $f(w_i)=f(w_j)$ for all $i,j\in [k]$, $f(y_i)\neq f(y_j)$ for all $i\neq j\in [k]$ and $f(z_i)\neq f(z_j)$ for all $i\neq j\in [k]$. %Call such an $f$ fat if the codegree of the vertices $f(y_1),\dots,f(y_k),f(w_1)$ is at least $\log n$.
    Observe that the number of such homomorphisms in which the vertices $f(y_1),\dots,f(y_k),f(w_1)$ have codegree at most $\log n$ is at most $\hom^*(K_{k,1}\Box K_2,G)(\log n)^{k-1}$. Indeed, there are at most $\hom^*(K_{k,1}\Box K_2,G)$ ways to choose $f|_{\{x_1,\dots,x_k,y_1,\dots,y_k,w_1,z_1\}}$, and given such a choice with the vertices $f(y_1),\dots,f(y_k),f(w_1)$ having codegree at most $\log n$, the number of ways to choose $f(z_2),\dots,f(z_k)$ is at most $(\log n)^{k-1}$, noting that we require that $f(w_1)=\dots=f(w_k)$. By Lemma \ref{lem:bounds for homs} \ref{item: (a)}, we have $\hom^*(K_{k,1}\Box K_2,G)(\log n)^{k-1}<n^{-\eps/2}\hom^*(K_{k,k}\Box K_2,G)$, so an $n^{-o(1)}$ proportion of all $f\in \Hom^*(K_{k,k}\Box K_2,G)$ satisfies that $f(w_i)=f(w_j)$ for all $i,j\in [k]$, $f(y_i)\neq f(y_j)$ for all $i\neq j\in [k]$,  $f(z_i)\neq f(z_j)$ for all $i\neq j\in [k]$ and the codegree of the vertices $f(y_1),\dots,f(y_k),f(w_1)$ is at least $\log n$.
    By Lemma \ref{lem:hom inequalities}, the number of such homomorphisms is therefore at least $n^{-o(1)} h^{k^2} (nd)^{2k-2k^2}$. For any such homomorphism, the vertices $f(y_1),\dots,f(y_k),f(w_1),f(z_1),\dots,f(z_k)$ are distinct since $f(y_i)\neq f(y_j)$ and $f(z_i)\neq f(z_j)$ for all $i\neq j\in [k]$ and $f$ is weakly injective. Hence, these vertices form a fat $K_{k,k+1}$. The number of times we count each such fat $K_{k,k+1}$ is at most $O(\lambda^k)$ since the number of choices for each $f(x_i)$ is at most $\lambda$, as we are only counting light homomorphisms $f$. It follows that there are at least $n^{-o(1)} h^{k^2} (nd)^{2k-2k^2} / \lambda^k$ copies of fat $K_{k,k+1}$.
\end{proof}

\subsection{Embedding procedure 2} \label{sec:embedding blowup}

In this subsection we define our second embedding procedure and prove Proposition \ref{prop:embedding 2}.
We need the following lemma, which is inspired by a result in \cite{grzesik2022turan}.

\begin{lemma} \label{lemma:blowup from ktt}
    For any positive integers $k$ and $\ell$, there exists some $C$ such that the following holds. Let $G$ be a graph in which the number of copies of $K_{k+1,k+1}$ is greater than $C$ times the number of copies of $K_{1,k+1}$. Then $G$ contains the $(k+1)$-blowup of a path of length $\ell$.
\end{lemma}

\begin{proof}
    We will show that a sufficiently large $C$ works. Define an auxiliary graph $\mathcal{G}$ whose vertices are the $(k+1)$-element subsets of $V(G)$ with non-empty common neighbourhood and in which there is an edge between $S\in V(G)^{(k+1)}$ and $T\in V(G)^{(k+1)}$ if $S\cap T=\emptyset$ and there is a complete bipartite graph between $S$ and $T$. Notice that $|V(\mathcal{G})|$ is at most the number of copies of $K_{1,k+1}$ in $G$, and that $e(\mathcal{G})$ is equal to the number of copies of $K_{k+1,k+1}$ in $G$. Therefore, by assumption, we have $e(\mathcal{G})>C|V(\mathcal{G})|$.

    Consider a random walk $(S_i)_{0\leq i\leq \ell}\in V(\mathcal{G})$, where $S_0$ is chosen randomly according to the stationary distribution of $\mathcal{G}$ and $S_{i+1}$ is a uniformly random neighbour of $S_i$. Note that each $S_i$ follows the stationary distribution. Call a vertex $T\in V(\mathcal{G})$ \emph{poor} if its degree in $\mathcal{G}$ is at most $\frac{C}{(\ell+1)^2}$. Now note that for each $0\leq i\leq \ell$, the probability that $S_i$ is poor is
    $$\sum_{T \text{ is poor}} d_{\mathcal{G}}(T)/(2e(\mathcal{G}))\leq |V(\mathcal{G})|\frac{C}{(\ell+1)^2}/(2e(\mathcal{G}))<\frac{1}{2(\ell+1)^2}.$$

    Now fix some $1\leq j\leq \ell$. Let $\mathcal{E}$ be the event that there is some $0\leq i\leq j-1$ such that $S_i$ is poor. By the above, we have $\mathbb{P}(\mathcal{E})< \frac{1}{2(\ell+1)}$. Note that
    \begin{equation}
        \mathbb{P}(S_j\cap \bigcup_{0\leq i\leq j-1} S_i\neq \emptyset)=\mathbb{P}(S_j\cap \bigcup_{0\leq i\leq j-1} S_i\neq \emptyset \wedge \mathcal{E})+\mathbb{P}(S_j\cap \bigcup_{0\leq i\leq j-1} S_i\neq \emptyset \wedge \bar{\mathcal{E}}). \label{eqn:random walk two terms}
    \end{equation}
    Observe that $\mathbb{P}(S_j\cap \bigcup_{0\leq i\leq j-1} S_i\neq \emptyset \wedge \mathcal{E})\leq \mathbb{P}(\mathcal{E})< \frac{1}{2(\ell+1)}$. We will show that also $\mathbb{P}(S_j\cap \bigcup_{0\leq i\leq j-1} S_i\neq \emptyset \wedge \bar{\mathcal{E}})<\frac{1}{2(\ell+1)}$. It suffices to prove that conditioning on any outcomes of $S_0,\dots,S_{j-1}$ such that $\mathcal{E}$ does not hold, the probability that $S_j\cap \bigcup_{0\leq i\leq j-1} S_i\neq \emptyset$ is less than $\frac{1}{2(\ell+1)}$. Note that $S_{j}$ is a random neighbour of $S_{j-1}$, that is, a random $(k+1)$-element subset of the common neighbourhood of $S_{j-1}$ in $G$, which has size at least $(C/(\ell+1)^2)^{1/(k+1)}$ since $S_{j-1}$ is not poor. Since $C$ is sufficiently large and $|\bigcup_{0\leq i\leq j-1} S_i|\leq \ell(k+1)$, it follows that indeed the conditional probability that $S_j\cap \bigcup_{0\leq i\leq j-1} S_i\neq \emptyset$ is less than $\frac{1}{2(\ell+1)}$. Thus, by (\ref{eqn:random walk two terms}), we have $\mathbb{P}(S_j\cap \bigcup_{0\leq i\leq j-1} S_i\neq \emptyset)<\frac{1}{\ell+1}$. Taking a union bound over all $j$, we conclude that with positive probability the sets $S_0,\dots,S_{\ell}$ are pairwise disjoint, in which case they span a $(k+1)$-blowup of a path of length $\ell$, as desired.
\end{proof}

We are now ready to prove Proposition \ref{prop:embedding 2}.

\begin{proof}[Proof of Proposition \ref{prop:embedding 2}]
We will show that a sufficiently large $C$ works. Let $N_i$ be the number of copies of $K_{i,k+1}$ in $G$. By assumption, $N_j>C N_1$. Let $N_i'$ be the number of copies of $K_{i,k+1}$ in which the part of size $k+1$ has codegree at most $k$. Clearly, $N_j'\leq N_1' k^{j-1}\leq N_1 k^{j-1}$. It follows that $N_j-N_j'>CN_1/2$. On the other hand, clearly, $N_{k+1}>c_k(N_j-N_j')$ for some positive $c_k$ depending only on $k$, since each $K_{j,k+1}$ in which the part of size $k+1$ has codegree greater than $k$ extends to a copy of $K_{k+1,k+1}$. Since $C$ is sufficiently large, it follows that $N_{k+1}>C'N_1$ for some sufficiently large $C'$. Now an application of Lemma \ref{lemma:blowup from ktt} provides the desired conclusion.
\end{proof}

\subsection{Embedding procedure 3}

In this subsection, we define our third embedding procedure and prove Proposition \ref{prop:embedding 3}. The following definition will be convenient to use.

\begin{definition}
    A copy of $K_{k,k+1}$ in a graph $G$ is called $L$-controlled if the codegree of the part of size $k+1$ is at least $L$ and less than $2L$. 
\end{definition}

We now prove the following ``cleaning lemma''.

\begin{lemma} \label{lem:cleaning for chain}
    Let $L$, $q$ and $k$ be positive integers satisfying $L\geq 4(k+1)q$. Let $G$ be a bipartite graph with parts $A$ and $B$ such that the number of $L$-controlled copies of $K_{k,k+1}$ with the part of size $k$ in $A$ is greater than $10(k+1)q$ times the number of copies of $K_{k,k}$. Then there exists a non-empty collection $\mathcal{F}$ of copies of $K_{k+1,k+1}$ in $G$ such that for any $F\in \mathcal{F}$ and any $k$-element subsets $S\subset V(F)\cap A$ and $T\subset V(F)\cap B$, there exist distinct vertices $a_1,\dots,a_q\in A\setminus S$ and $b_1,\dots,b_q\in B\setminus T$ such that $G[S\cup \{a_i\},T\cup \{b_i\}]\in \mathcal{F}$ for all $i\in [q]$.
\end{lemma}

\begin{proof}
    Call a copy of $K_{k+1,k+1}$ in $G$ $L$-controlled if its part in $B$ has codegree at least $L$ and less than $2L$. 
    
    Let us write $N_{k,k}$ for the number of copies of $K_{k,k}$, let $N_{k,k+1}$ be the number of $L$-controlled copies of $K_{k,k+1}$ with the part of size $k$ in $A$ and let $N_{k+1,k+1}$ be the number of $L$-controlled copies of $K_{k+1,k+1}$. By the assumption of the lemma, we have $N_{k,k+1}> 10(k+1)q N_{k,k}$.

    Notice that every $L$-controlled copy of $K_{k,k+1}$ with the part of size $k$ in $A$ extends to at least $L-k\geq L/2$ many $L$-controlled copies of $K_{k+1,k+1}$.
    Since each $L$-controlled copy of $K_{k+1,k+1}$ is obtained from $k+1$ different $L$-controlled copies of $K_{k,k+1}$, it follows that $N_{k+1,k+1}\geq \frac{L}{2(k+1)}N_{k,k+1}$. In particular, we have $N_{k+1,k+1}>5qL N_{k,k}$.

    Let $\mathcal{F}_0$ be the collection of $L$-controlled copies of $K_{k+1,k+1}$ in $G$. Having defined $\mathcal{F}_i$, we either terminate the process and set $\mathcal{F}=\mathcal{F}_i$, or define $\mathcal{F}_{i+1}\subsetneq \mathcal{F}_i$, according to the following rules.

    \begin{enumerate}
        \item [(a)] If there exists some $F\in \mathcal{F}_i$ such that for some $k$-element subsets $S\subset V(F)\cap A$ and $T\subset V(F)\cap B$, there are at most $2qL$ elements of $\mathcal{F}_i$ extending $S\cup T$, then let $\mathcal{F}_{i+1}$ be obtained from $\mathcal{F}_i$ by discarding all elements of $\mathcal{F}_i$ containing $S\cup T$.

        \item [(b)] If (a) does not hold, but there exists some $F\in \mathcal{F}_i$ such that for some $k$-element subset $S\subset V(F)\cap A$ and for $T= V(F)\cap B$, there are at most $\frac{L}{4(k+1)}$ elements of $\mathcal{F}_i$ extending $S\cup T$, then let $\mathcal{F}_{i+1}$ be obtained from $\mathcal{F}_i$ by discarding all elements of $\mathcal{F}_i$ containing $S\cup T$.
    \end{enumerate}

    If neither (a) nor (b) holds, then we terminate the process and set $\mathcal{F}=\mathcal{F}_i$. Since we discard at least one element in each step, the process eventually terminates.

    \medskip

    \noindent \textit{Claim.} We have $\mathcal{F}\neq \emptyset$.

    \medskip

    \noindent \textit{Proof of Claim.} We bound the number of elements discarded by each of steps (a) and (b). Note that each time we perform deletion step (a), we discard at most $2qL$ elements. Moreover, the number of times we can perform such a step is at most $N_{k,k}$ since $G[S\cup T]$ forms a copy of $K_{k,k}$ in any such step. Hence, the number of elements discarded by step (a) throughout the process is at most $2qLN_{k,k}<N_{k+1,k+1}/2=|\mathcal{F}_0|/2$. Similarly, each time we perform deletion step (b), we discard at most $\frac{L}{4(k+1)}$ elements. Moreover, the number of times we can perform such a step is at most $N_{k,k+1}$ since in any such step $G[S\cup T]$ forms an $L$-controlled copy of $K_{k,k+1}$ with the part of size $k$ in $A$. Hence, the number of elements discarded by step (b) throughout the process is at most $\frac{L}{4(k+1)}N_{k,k+1}\leq N_{k+1,k+1}/2=|\mathcal{F}_0|/2$. It follows that the total number of elements discarded throughout the process is less than $|\mathcal{F}_0|$, so $\mathcal{F}$ is indeed non-empty. $\Box$

    \medskip

    We shall now prove that $\mathcal{F}$ satisfies the conclusion of the lemma. Indeed, let $F\in \mathcal{F}$ and let $S\subset V(F)\cap A$ and $T\subset V(F)\cap B$ be $k$-element subsets. Since the process terminated with $\mathcal{F}$, by (a) it follows that there are at least $2qL$ elements of $\mathcal{F}$ extending $S\cup T$. Each extension involves precisely one vertex in $A\setminus S$ and one vertex in $B\setminus T$. Note that there are at most $2L$ extensions involving the same vertex $v\in B\setminus T$ because for any such $v$, the codegree of the set $T\cup \{v\}$ is at most $2L$ (else, there is no $L$-controlled copy of $K_{k+1,k+1}$ containing $T\cup \{v\}$). Thus, there must exist distinct vertices $b_1,\dots,b_q\in B\setminus T$ appearing in at least one extension of $S\cup T$ to an element of $\mathcal{F}$. Since the process terminated with $\mathcal{F}$, by (b), for each $i\in [q]$ there are at least $\frac{L}{4(k+1)}\geq q$ elements of $\mathcal{F}$ extending $S\cup (T\cup \{b_i\})$. Thus, we may choose $a_1,\dots,a_q$ greedily.
\end{proof}

The following corollary of Lemma \ref{lem:cleaning for chain} is the heart of this embedding procedure. Recall that $F_{k,\ell}$ is the graph with vertex set $\{x_1,\dots,x_{\ell},y_1,\dots,y_{\ell}\}$ and edge set $\{x_iy_j: |i-j|\leq k\}$.

\begin{lemma} \label{lem:build chain}
    Let $L$, $\ell$ and $k$ be positive integers satisfying $L\geq 8(k+1)\ell$ and $\ell\geq k+1$. Let $G$ be a bipartite graph with parts $A$ and $B$ such that the number of $L$-controlled copies of $K_{k,k+1}$ with the part of size $k$ in $A$ is greater than $20(k+1)\ell$ times the number of copies of $K_{k,k}$. Then $G$ contains $F_{k,\ell}$.
\end{lemma}

\begin{proof}
    By Lemma \ref{lem:cleaning for chain} (applied with $q=2\ell$), there exists a non-empty collection $\mathcal{F}$ of copies of $K_{k+1,k+1}$ in $G$ such that for any $F\in \mathcal{F}$ and any $k$-element subsets $S\subset V(F)\cap A$ and $T\subset V(F)\cap B$, there exist distinct vertices $a_1,\dots,a_{2\ell}\in A\setminus S$ and $b_1,\dots,b_{2\ell}\in B\setminus T$ such that $G[S\cup \{a_i\},T\cup \{b_i\}]\in \mathcal{F}$ for all $i\in [2\ell]$. We build a copy of $F_{k,\ell}$ in $G$ as follows. Choose vertices $x_1,\dots,x_{k+1}$ and $y_1,\dots,y_{k+1}$ forming a $K_{k+1,k+1}$ in $\mathcal{F}$. Now assume that we have already defined $x_1,\dots,x_t\in A$ and $y_1,\dots,y_t\in B$ suitably for some $k+1\leq t\leq \ell-1$ and that the $K_{k+1,k+1}$ formed by the vertices $x_{t-k},x_{t-k+1},\dots,x_t$ and $y_{t-k},y_{t-k+1},\dots,y_t$ belongs to $\mathcal{F}$. Let $S=\{x_{t-k+1},x_{t-k+2},\dots,x_t\}$ and let $T=\{y_{t-k+1},y_{t-k+2},\dots,y_t\}$. By the defining property of $\mathcal{F}$, there exist distinct vertices $a_1,\dots,a_{2\ell}\in A\setminus S$ and $b_1,\dots,b_{2\ell}\in B\setminus T$ such that $G[S\cup \{a_i\},T\cup \{b_i\}]\in \mathcal{F}$ for all $i\in [2\ell]$. In particular, there exist some $x_{t+1}\in A\setminus \{x_1,\dots,x_{t}\}$ and $y_{t+1}\in B\setminus \{y_1,\dots,y_{t}\}$ such that $G[\{x_{t-k+1},\dots,x_{t+1}\},\{y_{t-k+1},\dots,y_{t+1}\}]\in \mathcal{F}$. Hence, we can iterate this process and build a copy of $F_{k,\ell}$ in $G$.
\end{proof}

We are now ready to prove Proposition \ref{prop:embedding 3}.

\begin{proof}[Proof of Proposition \ref{prop:embedding 3}]
    We will show that it suffices to prove the result for $\ell\ge k+1$ with $C=80(k+1)\ell$. Indeed, if $\ell\le k$, then we apply the result with the length parameter equal to $k+1$, obtaining a copy of $F_{k,k+1}$, and then restrict it to $F_{k,\ell}$; this allows us to take $C=80(k+1)^2$. Thus, we may assume that $\ell\ge k+1$. Let $A$ and $B$ be the parts of $G$. Without loss of generality, the number of fat copies of $K_{k,k+1}$ with the part of size $k$ in $A$ is greater than $\frac{C}{2}\log n$ times the number of copies of $K_{k,k}$. Note that any such copy is $L$-controlled for some $L=2^i(\log n)$ and nonnegative integer $i\leq \log_2 n$. It follows that for some $L$ of this form, the number of $L$-controlled copies of $K_{k,k+1}$ with the part of size $k$ in $A$ is greater than $C/4$ times the number of copies of $K_{k,k}$. Since $C/4=20(k+1)\ell$, $L\geq \log n$ and $n$ is sufficiently large, Lemma \ref{lem:build chain} implies that $G$ contains $F_{k,\ell}$.
\end{proof}

\section{Concluding remarks}

The main remaining problem is to extend Theorem \ref{thm:even min degree} to odd values of $s$.  An even stronger direction is to ask for a small \emph{$s$-regular} subgraph.  The optimal bound has already been proved in the bounded-order setting for even $s$ \cite{janzer2023rainbow} and for $s=3$ \cite{janzer2023disproof}. This problem can be viewed as an interpolation between two important topics: the Tur\'an number of $s$-regular graphs (where a bounded-size regular subgraph is sought) and the Erd\H{o}s--Sauer problem~\cite{janzer2023regular} (where we are satisfied with any regular subgraph). 

As a less ambitious goal, it would already be interesting to prove the strengthening of Theorem~\ref{thm:small odd subgraph} where a bounded-order $s$-regular subgraph is found. Finally, it would also be interesting to replace the family $\mathcal{H}$ in Theorem \ref{thm:odd extremal} by a single graph.

\medskip

\paragraph{Statement on AI use.} The results proved in this paper were obtained without AI assistance. We used ChatGPT Pro to polish an earlier draft.

\bibliographystyle{abbrv}
\bibliography{bibliography}

\end{document}